\documentclass{elsarticle}

\usepackage{lineno,hyperref}
\modulolinenumbers[5]

\journal{}
\usepackage{geometry}
\usepackage{amsmath}
\usepackage{subfigure}
\usepackage{latexsym}
\usepackage{amsmath}
\usepackage{amssymb}
\usepackage{amsfonts}
\usepackage{verbatim}
\usepackage[latin1]{inputenc}
\usepackage{mathrsfs}
\usepackage{amsfonts}
\usepackage{graphicx}
\usepackage{colortbl,dcolumn}
\usepackage{amsmath}
\usepackage{psfrag}
\usepackage{lipsum}
\usepackage{amsfonts}
\usepackage{graphicx}
\usepackage{epstopdf}
\usepackage{algorithmic}
\usepackage{ulem}
\ifpdf
  \DeclareGraphicsExtensions{.eps,.pdf,.png,.jpg}
\else
  \DeclareGraphicsExtensions{.eps}
\fi
\usepackage[misc]{ifsym}
\usepackage{enumerate}

\newtheorem{ap}{Assumption}
\newtheorem{tm}{Theorem}
\newtheorem{rk}{Remark}
\newtheorem{df}{Definition}
\newtheorem{prop}{Proposition}
\newtheorem{lm}{Lemma}
\newtheorem{cor}{Corollary}

\newcommand{\E}{\mathbb E}

\newcommand{\N}{\mathbb N}

\newcommand{\bi}{\mathbf i}
\newcommand{\bs}{\mathbf s}

\newcommand{\OO}{\mathcal O}
\newcommand{\HH}{\mathbb H}

\newcommand{\<}{\langle}
\renewcommand{\>}{\rangle}
\usepackage{graphicx}
\allowdisplaybreaks[4]

\begin{document}

\begin{frontmatter}
	
	\title{Stochastic logarithmic Schr\"odinger equations: energy regularized approach}
	\tnotetext[mytitlenote]{This work was funded by National Natural Science Foundation of China (No.
		91630312, No. 91530118, No.11021101 and No. 11290142).}
	
	\author[c]{Jianbo Cui}
	\ead{jianbocui@lsec.cc.ac.cn}

	\author[cas]{Liying Sun\corref{cor}}
	\ead{liyingsun@lsec.cc.ac.cn}
	
	\cortext[cor]{Corresponding author.}
	
	\address[c]{School of Mathematics, Georgia Tech, Atlanta, GA 30332, USA}
	
	\address[cas]{1. LSEC, ICMSEC, 
		Academy of Mathematics and Systems Science, Chinese Academy of Sciences, Beijing,  100190, China\\
		2. School of Mathematical Science, University of Chinese Academy of Sciences, Beijing, 100049, China}

	\begin{abstract}
	In this paper, we prove the global existence and uniqueness of the solution of the stochastic logarithmic Schr\"odinger (SlogS) equation driven by additive noise or multiplicative noise. The key ingredient lies on the regularized stochastic logarithmic Schr\"odinger (RSlogS) equation with regularized energy and the 
	strong convergence analysis of the solutions of  (RSlogS) equations. 
	In addition, temporal H\"older regularity estimates and uniform estimates in energy space $\mathbb H^1(\mathcal O)$ and weighted Sobolev space $L^2_{\alpha}(\mathcal O)$ of the solutions for both SlogS equation and RSlogS equation are also obtained.
	\end{abstract}
	
	
\end{frontmatter}


	\section{Introduction}
	The deterministic logarithmic Schr\"odinger equation  
	has wide applications in quantum mechanics, quantum optics, nuclear physics, transport and diffusion phenomena, open quantum system, Bose--Einstein condensations and so on (see e.g. \cite{AZ03,BM76,Hef85,Lau08,Yas78,Zlo10}). 
		It takes the form of
		\begin{align*}
		\partial_t u(t,x)&=\mathbf i \Delta u(t,x)+\mathbf i \lambda u(t,x)\log(|u(t,x)|^2)+\bi V(t,x,|u|^2)u(t,x),\;  x\in \mathbb R^d,\; t>0,\\\nonumber
		u(0,x)&=u_0(x),\; \quad x\in \mathbb R^d,
		\end{align*} 
		where $\Delta$ is the Laplacian operator on $\mathcal{O} \subset \mathbb{R}^{d}$ with $\mathcal{O}$ being either $\mathbb{R}^{d}$ or a bounded domain with homogeneous Dirichlet or periodic boundary condition,
	$t$ is time, $x$ is spatial coordinate, $\lambda\in \mathbb R/\{0\}$ characterizes the force of nonlinear interaction, and $V$ is a real-valued function. While retaining many of the known features of the linear Schr\"odinger equation, Bialynicki--Birula and Mycielski show that only such a logarithmic nonlinearity satisfies the condition of
	separability of noninteracting systems (see \cite{BM76}). The logarithmic nonlinearity makes the logarithmic Schr\"odinger equation unique among nonlinear wave equations. 
	For instance, the longtime dynamics of the logarithmic Schr\"odinger equation is essentially different from the Schr\"odinger equation. There is a faster dispersive phenomenon when $\lambda<0$ and the convergence of the modulus of the solution to a universal Gaussian profile (see \cite{CG18}), and no dispersive phenomenon when $\lambda>0$ (see \cite{Caz83}).
	
	In this paper, we are mainly focus on the well-posedness of the following stochastic logarithmic Schr\"odinger (SlogS) equation,
		\begin{align}\label{SlogS}
		d u(t)&=\mathbf i \Delta u(t)dt+\mathbf i \lambda u(t)\log(|u(t)|^2)dt+ \widetilde g(u)\star dW(t),\; t>0\\\nonumber 
		u(0)&=u_0,
		\end{align}
		where  $W(t)=\sum_{k\in \mathbb N^+}Q^{\frac 12}e_k\beta_k(t)$, $\{e_k\}_{k\in\mathbb N^+}$ is an orthonormal basis of $L^2(\mathcal O;\mathbb C)$ with $\{\beta_k\}_{k\in\mathbb N^+}$ being a sequence of independent  Brownian motions on a probability space $(\Omega,$ 
		$\mathcal F, (\mathcal F_t)_{t\ge 0},\mathbb P).$ 
		Here $\widetilde g$ is a continuous function and $\widetilde g(u)\star dW(t)$ is defined by 
		\begin{align*}
		&\widetilde g(u)\star dW(t) =-\frac 12\sum_{k\in\mathbb N^+}|Q^{\frac 12}e_k|^2\Big(|g(|u|^2)|^2u\Big)dt\\\nonumber 
		&\quad-\bi \sum_{k\in \mathbb N^+}g(|u|^2)g'(|u|^2) |u|^2 u
		Im(Q^{\frac 12}e_k)Q^{\frac 12}e_k dt+ \mathbf i g(|u|^2)u dW(t)
		\end{align*}
		if $\widetilde g(x)=\bi g(|x|^2)x$, and by 
		$$\widetilde g(u)\star dW(t)= dW(t)$$ if $\widetilde g=1$. 
		We would like to remark that when $W$ is $L^2(\mathcal O;\mathbb R)$-valued and $\widetilde g(x)=\bi g(|x|^2)x$, $\widetilde g(u)\star dW(t)$ is just the classical Stratonovich integral. 
	
	The SlogS equation \eqref{SlogS} could be derived from the deterministic model by using Nelson's mechanics {\cite{Nel66}}. Applying the Madlung transformation $u(t,x)=\sqrt{\rho(t,x)}e^{\bi  {S(t,x)} }$, \cite{Nas86} obtains a fluid expression of the solution as follows,
	\begin{align*}
	\partial_tS(t,x)&=-|\nabla S(t,x)|^2- \frac {1}4 \frac {\delta I}{\delta \rho}(\rho(t,x)){+\lambda}\log(\rho){+V(t,x,\rho(t,x))},\; \\\nonumber
	\partial_t\rho(t,x)&=-2div(\rho(t,x) \nabla S(t,x)),\; S(0,x_0)=S_0(x_0)\;, \rho(0)=\rho_0,
	\end{align*}
	where $I(\rho)=\int_{\mathcal O}|\nabla \log(\rho)|^2$$\rho dx$ is the Fisher information. 
	If $V$ is random and fluctuates rapidly, the term $\mathbf iV u$ can be approximated by some multiplicative Gaussian noise $\widetilde g(u)\dot W,$ which plays an important role in the theory of measurements continuous in time in open quantum systems (see e.g. \cite{BG09}). Then we could use the inverse of Madlung transformation and formally obtain the stochastic logarithmic Schr\"odinger equation
	\begin{align*}
	\partial_t u(t,x)&=\mathbf i \Delta u(t,x)+\mathbf i \lambda u(t,x)\log(|u(t,x)|^2)+\widetilde g(u(t,x))\dot W(t,x),\; x\in \mathcal O,\; t>0\\
	u(0,x)&=u_0(x), \quad x\in \overline{\mathcal O}.
	\end{align*}
		The main assumption on $W$ and $\widetilde g$ is stated as follows.
		
		\begin{ap}\label{main-as}
			The diffusion operator is the Nemystkii operator of $\widetilde g$. Wiener process $W$ and $\widetilde g$ satisfies one of the following condition,
			\begin{enumerate}
				\item[{\rm Case 1.}] $\{W(t)\}_{t\ge 0}$ is $L^2(\mathcal O;\mathbb C)$-valued and $\widetilde g=1;$
				\item[{\rm Case 2.}] $\{W(t)\}_{t\ge 0}$ is $L^2(\mathcal O;\mathbb C)$-valued,  
				$\widetilde g(x)=\mathbf ig(|x|^2)x,$ $g \in \mathcal C^2_b(\mathbb R)$ and satisfies the growth condition 
				\begin{align*}
				\sup_{x\in [0,\infty)}|g(x)|+\sup_{x\in [0,\infty)}|g'(x)x|+\sup_{x\in [0,\infty)}|g''(x)x^2|\le C_g;
				\end{align*} 
				\item[{\rm Case 3.}] $ {\{W(t)\}_{t\ge 0}}$ is $L^2(\mathcal O;\mathbb R)$-valued, $\widetilde g(x)=\mathbf ig(|x|^2)x,$ $g \in \mathcal C^1_b(\mathbb R)$ and satisfies the growth condition 
				\begin{align*}
				\sup_{x\in [0,\infty)}|g(x)|+\sup_{x\in [0,\infty)}|g'(x)x|\le C_g.
				\end{align*} 
			\end{enumerate}
		\end{ap}

		\begin{ap}\label{ass-strong}
			Assume that $g$ satisfies
			\begin{align}\label{con-g}
			(x+y)(g(|x|^2)-g(|y|^2))\le C_g|x-y|, \quad x,y\in [0,\infty).
			\end{align}
			When $W(t)$ is $L^2(\mathcal O;\mathbb C)$-valued, we in addition assume that $g$ satisfies  following one-side Lipschitz continuity
			\begin{align}\label{con-g1}
			|(\bar y-\bar x)(g'(|x|^2)g(|x|^2)|x|^2x-g'(|y|^2)g(|y|^2)|y|^2y)| \le C_g|x-y|^2, \quad x,y\in \mathbb C.
			\end{align}
		\end{ap}
		A typical example is $\widetilde g(u)=u$, and then Eq. \eqref{SlogS}  becomes the SlogS equation driven by linear multiplicative noise in \cite{BRZ17}.
	
	There are two main difficulties in proving the well-posedness of the SlogS equation. 
	 On one hand, the random perturbation in SlogS equation destroys a lot of physical  conservation laws, like the mass and energy conservation laws in Case 1 and Case 2, and the energy conservation law in Case 3. 
	Similar phenomenon has been observed in stochastic nonlinear Schr\"odinger equation with polynomial nonlinearity (see \cite{BD03}). 
	 On the other hand, the logarithmic nonlinearity in SlogS equation is not locally Lipschitz  continuous. 
	The contraction mapping arguments via Strichartz estimates (see e.g. \cite{BRZ16,BD03,Hor18}) for stochastic nonlinear Schr\"odinger equation with smooth nonlinearity are not applicable here. 
	We only realize that in Case 2, when the driving noise is a linear multiplicative noise ($\widetilde g(u)=\bi u$), \cite{BRZ17} uses a rescaling technique, together with maximal monotone operator theory to obtain a unique global mild solution in some Orlicz space. 
	As far as we know, there has no results concerning the well-posedness of the SlogS equation driven by additive noise or general multiplicative noise.

	To show the well-posedness of the considered model, we introduce an energy regularized problem inspired by \cite{BCST19} where the authors use the regularized problem to study error estimates of numerical methods for deterministic logarithmic Schr\"odinger equation. 
	The main idea is firstly constructing a proper approximation of $\log(|x|^2)$ denoted by $f_{\epsilon}(|x|^2).$ Then it induces the regularized entropy 
	$F_{\epsilon}$ which is an approximation of the entropy $F(\rho)=\int_{\mathcal O}(\rho\log(\rho)-\rho)dx,$ where $\rho=|u|^2$. 
The RSlogS equation is defined by
		\begin{align}\label{Reg-SlogS}
		du^{\epsilon}=\mathbf i \Delta u^{\epsilon}dt+\mathbf i u^{\epsilon} f_{\epsilon}(|u^{\epsilon}|^2)dt+ \widetilde g(u^{\epsilon})\star dW(t),
		\end{align}
		whose regularized energy is $\frac 12\|\nabla u\|^2+\frac {\lambda}2 F_{\epsilon}(u)$. 
		Denoting $\mathbb H:=L^2=L^2(\mathcal O;\mathbb C)$ with the product $\<u,v\>:=\int_{\mathbb R^d} Re(u\bar v)dx,$ for $u,v\in \mathbb H,$
		we obtain the existence and uniqueness of the solution of regularized SlogS equation by proving $\epsilon$-independent estimate in $\mathbb H^1:=W^{1,2}$ and the weighted $L^2$-space $L^2_{\alpha}:=\{v\in L^2 | \; x \mapsto (1+|x|^2)^{\frac \alpha 2}v(x)\in L^2\}$ with  the norm $\|v\|_{L^2_\alpha(\mathbb R^d)}:=\|(1+|x|^2)^{\frac \alpha 2}v(x)\|_{L^2(\mathbb R^d)}$. 
		Then we are able to prove that the limit of $\{u^{\epsilon}\}_{\epsilon>0}$ is convergent to a unique stochastic process $u$ which is shown to be the unique mild solution of \eqref{SlogS}. Meanwhile, the sharp convergence rate of $\{u^{\epsilon}\}_{\epsilon>0}$ is given when $\mathcal O=\mathbb R^d,$ or $\mathcal O$ is a bounded domain in $\mathbb R^d$ equipped with homogenous Dirichlet or periodic boundary condition. 
		Our main result is formulated as follows.
		
		\begin{tm}\label{mild-general}
			Let $T>0$, Assumptions \ref{main-as} and  \ref{ass-strong} hold, $u_0\in \mathbb H^1\cap L^2_{\alpha},$ $\alpha\in(0,1],$ be $\mathcal F_0$ measurable and has any finite $p$th moment.   Assume that $\sum_{i\in\mathbb N^+} \|Q^{\frac 12}e_i\|_{L^2_{\alpha}}^2+\|Q^{\frac 12}e_i\|_{\mathbb H^1}^2 <\infty$ when $\widetilde g=1$ and that  $\sum_{i\in\mathbb N^+} \|Q^{\frac 12}e_i\|_{\mathbb H^1}^2+\|Q^{\frac 12}e_i\|_{W^{1,\infty}}^2<\infty$ when $\widetilde g(x)= \mathbf i g(|x|^2)x$.
			Then there exists a unique mild solution $u\in C([0,T];\mathbb H)$ of Eq. \eqref{SlogS}.
			Moreover,  for $p\ge 2$, there exists $C(Q,T,\lambda,p,u_0)>0$ such that 
			\begin{align*}
			\E\Big[\sup_{t\in [0,T]}\|u(t)\|_{\mathbb H^1}^p\Big]+\E\Big[\sup_{t\in [0,T]}\|u(t)\|_{L^2_{\alpha}}^p\Big]\le C(Q,T,\lambda, p,u_0).
			\end{align*}
		\end{tm}
		When $W(t)$ is $L^2(\mathcal O;\mathbb R)$-valued, the well-posedness of SlogS equation with a super-linearly growing diffusion coefficient is also proven (see Theorem \ref{mild-d1}).

	The reminder of this article is organized as follows. In section 2, we introduce  the RSlogS equation and show the local well-posedness of RSlogS equation driven by both additive and multiplicative noise. Section 3 is devoted to $\epsilon$-independent estimate of the mild solution in $\mathbb H^1$ and $L^2_{\alpha}$ of the RSlogS equation. In section 4, we prove the main result by passing the limit of the sequence of the regularized mild solutions and providing the sharp strong convergence rate. Several technique details are postponed to the Appendix. Throughout this article,  $C$ denotes
	various constants which may change from line to line.

	\section{Regularized SLogS equation}\label{sec-loc}

	In this section, we show the well-posedness of the solution for Eq. \eqref{Reg-SlogS} (see Appendix for the definition of the solution). We would like to remark that there are several choices of the regularization function $f_{\epsilon}(|x|^2)$. For instance, one may take $f_{\epsilon}(|x|^2)=\log(\frac{|x|^2+\epsilon}{1+\epsilon |x|^2})$ (see Lemma \ref{sec-reg} in Appendix for the necessary properties) or $f_{\epsilon}(|x|^2)=\log(|x|^2+\epsilon)$ (see e.g. \cite{BCST19} and references therein for more choices of regularization functions). If the regularization function $f_{\epsilon}$ enjoys the same properties of $\log(\frac{|x|^2+\epsilon}{1+\epsilon |x|^2})$, then one can follow our approach to obtain the well-posedness of Eq. \eqref{Reg-SlogS}. In the following, we first present the local well-posedness of Eq. \eqref{Reg-SlogS}, and then derive global existence and uniform estimate of its solution.  For simplicity, we always assume that $0<\epsilon \ll 1$.

	\subsection{Local well-posedness of Regularized SLogS equation}
	
	In this part, we give the detailed estimates to get the local well-posedness in $\mathbb H^2$ of Eq. \eqref{Reg-SlogS}  if $d\le 3$ via the regularization function $f_{\epsilon}(|x|^2)=\log(|x|^2+\epsilon).$ 
	In the case of $d\ge 3$, one could use the regularization function like $f_{\epsilon}(|x|^2)=\log(\frac{|x|^2+\epsilon}{1+\epsilon |x|^2})$ to get the local well-posedness in $\mathbb H^1$.  Assume that $u_0\in \mathbb H^2$ when using the regularization $\log(|x|^2+\epsilon)$ and that $u_0\in \mathbb H^1$ when applying the regularization $\log(\frac{|x|^2+\epsilon}{1+\epsilon |x|^2}).$

	Denote by $\mathbb M_{\mathcal F}^p(\Omega;C([0,T];\mathbb H^2))$ with $p\in [1,\infty)$ the space of process $v: [0,T]\times \Omega \to \mathbb H^2$ with continuous paths in $\mathbb H^2$ which is $\mathcal F_t$-adapted and satisfies
	\begin{align*}
	\|v\|^p_{\mathbb M_{\mathcal F}^p(\Omega;C([0,T];\mathbb H^2))}
	:=\E\Big[\sup_{t\in[0,T]}\|u(t)\|_{\mathbb H^2}^p\Big]<\infty.
	\end{align*}
	Let $\tau\le T$ be an $\mathcal F_t$-stopping time.
	And we call  $v\in \mathbb M_{\mathcal F}^p(\Omega;C([0,\tau);\mathbb H^2))$, if there exists $\{\tau_n\}_{n\in \mathbb N^+}$ with $\tau_n\nearrow \tau$ as $n\to \infty$ a.s., such that $v\in \mathbb M_{\mathcal F}^p(\Omega;C([0,\tau_n];\mathbb H^2))$ for $n\in \mathbb N^+.$ 
	Next we show the existence and uniqueness of the local mild solution (see Definition \ref{global} in Appendix).
	
		For the sake of simplicity, let us ignore the dependence on $\epsilon$ and write $u_R:=u^{\epsilon}_R,$ where $u_R$ is the solution of the truncated equation 
		\begin{align}\label{trun-reg}
		du_R=&\mathbf i \Delta u_R dt+\mathbf i \lambda \Theta_R(u_R,t) u_Rf_{\epsilon}(|u_R|^2)dt\\\nonumber
		&-\frac 12 \Theta_R(u_R,t) \sum_{k\in\mathbb N^+}|Q^{\frac 12}e_k|^2\Big(|g(|u_R|^2)|^2u_R\Big)dt\\\nonumber
		& -\bi \Theta_R(u_R,t)\sum_{k\in \mathbb N^+}g(|u_R|^2)g'(|u_R|^2) |u_R|^2 u_R
		Im(Q^{\frac 12}e_k)Q^{\frac 12}e_kdt\\\nonumber
		&+ \Theta_R(u_R,t) \mathbf i g(|u_R|^2)u_R dW(t).
		\end{align}
		Here, $\Theta_R(u,t):=\theta_R(\|u\|_{\mathcal C([0,t];\mathbb H^2)}), R>0,$ with
		a cut-off function $\theta_R$, that is, a positive $C^{\infty}$ function on $\mathbb R^+$ which has a compact support, and 
		\begin{equation*}
		\theta_R(x)=
		\left\{
		\begin{aligned}
		&0, \quad \text{for}\quad x\ge 2R,\\
		&1,\quad \text{for}\quad x\in [0,1].
		\end{aligned}
		\right.
		\end{equation*}

	\begin{lm}\label{trun}
		Let  Assumption \ref{main-as} hold, $d\le 3$, and $f_{\epsilon}(|x|^2)=\log(|x|^2+\epsilon)$. Assume in addition that $g \in \mathcal C^3_b(\mathbb R)$ when $W(t)$ is $L^2(\OO;\mathbb C)$-valued and that $g \in \mathcal C^2_b(\mathbb R)$ when $W(t)$ is $L^2(\OO;\mathbb R)$-valued.  
		Suppose that the $Q$-Wiener process $W(t)$ satisfies $\sum_{i\in\mathbb N^+}\|Q^{\frac 12}e_i\|_{\mathbb H^2}^2<\infty$, and $u_0\in \mathbb H^2$ is $\mathcal F_0$-measurable and has any finite $p$th moment. Then there exists a unique global solution to \eqref{trun-reg} with continuous $\mathbb H^2$-valued path.
	\end{lm}

	{\noindent \bf{Proof}}
			Let $S(t)=\exp(\bi\Delta t)$ be the $C_0$-group generated by $\bi\Delta.$
		For fixed $R>0$,  we use the following notations, for $t\in [0,T],$
		{\small
			\begin{align*}
			\Gamma_{det}^R u(t):&=\mathbf i\int_0^t S(t-s) \Big( \Theta_R(u,s)\lambda f_{\epsilon}(|u(s)|^2)u(s)\Big)ds,\\
			\Gamma_{mod}^R u(t):&=-\frac 12\int_0^t S(t-s)\Big(\Theta_R(u,s)\sum_{k\in\mathbb N^+}|Q^{\frac 12}e_k|^2\Big(|g(|u(s)|^2)|^2u(s)
			\Big)\Big)ds\\
			-\bi \int_0^tS(t-s)&\Big(
			\Theta_R(u,s)\sum_{k\in \mathbb N^+}g(|u(s)|^2)g'(|u(s)|^2) |u(s)|^2 u(s)
			Im(Q^{\frac 12}e_k)Q^{\frac 12}e_k \Big)ds,\\
			\Gamma_{Sto}^R u(t):&=\mathbf i\int_0^t S(t-s) \Big(\Theta_R(u,s)g(|u(s)|^2)u(s)\Big)dW(s).
			\end{align*}
			W}e look for a fixed point of the following operator given by 
		\begin{align*}
		\Gamma^R u(t):=S(t)u_0+\Gamma_{det}^R u(t)+\Gamma_{mod}^R u(t)
		+\Gamma_{Sto}^R u(t), u\in \mathbb M_{\mathcal F}^p(\Omega; C([0,r];\mathbb H^2)),
		\end{align*}
		where $r$ will be chosen later. 
		The unitary property of $S(\cdot)$ yields that
		$$\|S(\cdot)u_0\|_{\mathbb M_{\mathcal F}^p(\Omega; C([0,r];\mathbb H^2))}\le \|u_0\|_{\mathbb H^2}.$$
		Now, we define a stopping time $\tau=\inf\{t\in [0,T]: \|u\|_{C([0,t];\mathbb H^2)}\ge 2R\}\wedge r.$
		By using the properties of $\log(\cdot)$ and the Sobolev embedding $\mathbb H^2\hookrightarrow L^\infty(\OO)$, we have 
		\begin{align*}
		\|\Gamma_{det}^R u\|_{C([0,r];\mathbb H^2)}&\le Cr|\lambda| \max(|\log(\epsilon)|,\log(\epsilon+4R^2))(\sup_{t\in [0,r]}\|u(t)\|_{\mathbb H^2})\\
		&\quad+Cr|\lambda|(1+\epsilon^{-1})\Big(\sup_{t\in [0,r]}\|u(t)\|_{\mathbb H^2}+\sup_{t\in [0,r]}\|u(t)\|_{\mathbb H^2}^{2}\Big)\\
		&\le C(\epsilon,\lambda)r(1+2R+R^2),
		\end{align*}
		and 
		\begin{align*}
		&\|\Gamma_{mod}^R u\|_{C([0,r];\mathbb H^2)}
		\le C(\lambda,g)r\sum_{k\in \N^+}\|Q^{\frac 12}e_k\|_{\mathbb H^2}^2\Big(R+R^3+R^5\Big).
		\end{align*}
		Integrating over $\Omega$ yields that 
		\begin{align*}
		\|\Gamma_{det}^R u\|_{\mathbb M_{\mathcal F}^p(\Omega;C([0,r];\mathbb H^2))}+\|\Gamma_{mod}^R u\|_{\mathbb M_{\mathcal F}^p(\Omega;C([0,r];\mathbb H^2))}
		\le& C(\lambda,g)r\sum_{k\in \N^+}\|Q^{\frac 12}e_k\|_{\mathbb H^2}^2\Big(R+R^3+R^5\Big).
		\end{align*} 
		The Burkerholder inequality yields that for $p\ge 2$,
		\begin{align*}
		\|\Gamma_{sto}^R u\|_{\mathbb M_{\mathcal F}^p(\Omega; C([0,r];\mathbb H^2))}
		\le& C(g) \Big(\E\Big[\Big(\int_0^r\sum_{k\in\mathbb N^+} \|g(|u(s)|^2)u(s)Q^{\frac 12}{e_k}\|_{\mathbb H^2}^2ds\Big)^{\frac p2}\Big]\Big)^{\frac 1p}\\
		\le& C(g)r^{\frac 12}\left(\sum_{k\in \N^+}\|Q^{\frac 12}e_k\|_{\mathbb H^2}^2\right)^{\frac 12}\Big(R+R^3\Big).
		\end{align*} 
		Therefore, $\Gamma^R$ is well-defined on $\mathbb M_{\mathcal F}^p(\Omega; C([0,r];\mathbb H^2)).$ 
		
		Now we turn to show the contractivity of $\Gamma^R.$
		Let $u_1, u_2\in \mathbb M_{\mathcal F}^p(\Omega;C([0,r];\mathbb H^2)),$ and define the stopping times
		$\tau_j=\inf\{t\in [0,T]: \|u_j\|_{C([0,r];\mathbb H^2)}\ge2R\}\wedge r, j=1,2.$
		For a fixed $\omega,$ let us assume that $\tau_1\le \tau_2$ without the loss of generality. 
		Then direct calculation leads to 
		\begin{align*}
		&\|\Gamma_{det}^R u_1-\Gamma_{det}^R u_2\|_{C([0,r];\mathbb H^2)}\\
		\le &
		C(\lambda)r \Big\|\Big(\Theta_R(u_2,\cdot)-\Theta_R(u_1,\cdot)\Big) f_{\epsilon}(|u_2|^2)u_2\Big\|_{C([0,r];\mathbb H^2)} \\
		&+C(\lambda)r  \Big\|\Theta_R(u_1,\cdot) \Big( f_{\epsilon}(|u_1|^2)u_1- f_{\epsilon}(|u_2|^2)u_2\Big)\Big\|_{C([0,r];\mathbb H^2)}\\
		\le& C(\lambda)r\|u_1-u_2\|_{C([0,r];\mathbb H^2)} \|f_{\epsilon}(|u_2|^2)u_2\|_{C([0,\tau_2];\mathbb H^2)}\\
		&\quad+C(\lambda)r\|f_{\epsilon}(|u_2|^2)u_2-f_{\epsilon}(|u_1|^2)u_1\|_{C([0,\tau_1];\mathbb H^2)}\\
		\le& C(\lambda,\epsilon)r\|u_1-u_2\|_{C([0,r];\mathbb H^2)}(1+R^2),
		\end{align*}
		and 
		\begin{align*}
		\|\Gamma_{mod}^R u_1-\Gamma_{mod}^R u_2\|_{C([0,r];\mathbb H^2)}
		\le C(\lambda,\epsilon) \sum_{k\in\mathbb N^+}\|Q^{\frac 12}e_k\|_{\mathbb H^2}^2 r\|u_1-u_2\|_{C([0,r];\mathbb H^2)}(1+R^6).
		\end{align*}
		By applying the Burkerholder inequality, we obtain
		\begin{align*}
		&\|\Gamma_{Sto}^R u_1-\Gamma_{Sto}^R u_2\|_{\mathbb M_{\mathcal F}^p(\Omega; C([0,r];\mathbb H^2))}\\
		\le& Cr \Big\|\Theta_R(u_1,\cdot)(g(|u_1|^2)u_1-g(|u_2|^2)u_2) \Big\|_{L^p(\Omega; L^2([0,r];\mathbb L_2^{2}))}\\
		&+  Cr \Big\|(\Theta_R(u_1,\cdot)-\Theta_R(u_2,\cdot))g(|u_2|^2)u_2 \Big\|_{L^p(\Omega; L^2([0,r];\mathbb L_2^{2}))}\\
		\le& C(\lambda,\epsilon) (\sum_{k\in\mathbb N^+}\|Q^{\frac 12}e_k\|_{\mathbb H^2}^2)^{\frac 12} r^{\frac 12}\|u_1-u_2\|_{\mathbb M_{\mathcal F}^p(\Omega; C([0,r];\mathbb H^2))}(1+R^3).
		\end{align*}
		{ where $\mathbb L_2^{2}$ is the space of Hilbert--Schmidt operators form $U_0=Q^\frac 12(\mathbb L^2(\OO))$ to $\mathbb H^2.$} 
		Combining all the above estimates, we have 
		\begin{align*}
		&\quad\|\Gamma^R u_1-\Gamma^R u_2\|_{\mathbb M_{\mathcal F}^p(\Omega; C([0,r];\mathbb H^2))}\\
		&\le C(\lambda,\epsilon)\left(1+\sum_{k\in\mathbb N^+}\|Q^{\frac 12}e_k\|_{\mathbb H^2}^2\right) (r^{\frac 12}+r)\|u_1-u_2\|_{\mathbb M_{\mathcal F}^p(\Omega; C([0,r];\mathbb H^2))}(1+R^6),
		\end{align*}
		which implies that there exists a small $r>0$ depending on $Q,R,\lambda,\epsilon$ such that 
		$\Gamma^R$ is a strict contraction in $\mathbb M_{\mathcal F}^p(\Omega; C([0,r];\mathbb H^2))$ and has a fixed point $u^{R,1}$ satisfying $\Gamma^R(u^{R,1})=u^{R,1}.$
		
		Assume that we have found the fixed point  on each interval $[(l-1)r,l r],$ $l\le k$ for some $k\ge 1$. 
		Define
		\begin{align*}  
		u^{R,k}=S(\cdot )u_0+ \Gamma^{R}_{det}u^{R,k}+\Gamma^{R}_{mod}u^{R,k}+\Gamma^{sto}_{det}u^{R,k},
		\end{align*}
		on $[0,kr].$
		In order to extend $u^{R,k}$ to $[kr,(k+1)r]$, we repeat the previous arguments to show that on the interval $[kr,(k+1)r],$ there exists a fixed point of the map $\Gamma^R$ defined by
		\begin{align*}
		\Gamma^R u(t):=S(t)u^{R,k}({kr})+\Gamma_{det}^{R,k} u(t)+\Gamma_{mod}^{R,k} u(t)
		+\Gamma_{Sto}^{R,k} u(t), u\in \mathbb M_{\mathcal F}^p(\Omega; C([0,r];\mathbb H^2)).
		\end{align*}
		Here we use the following notations, 
		{\small 
			\begin{align*}
			\Gamma_{det}^{R,k}  u(t):&=\mathbf i\int_0^t S(t-s) \Big( \widetilde \Theta_R(u,k,s)\lambda f_{\epsilon}(|u(s)|^2)u(s)\Big)ds,\\
			\Gamma_{mod}^{R,k}  u(t):&=-\frac 12\int_0^t S(t-s)\Big(\widetilde \Theta_R(u,k,s)\sum_{j\in\mathbb N^+}|Q^{\frac 12}e_j|^2\Big(|g(|u(s)|^2)|^2u(s)\Big)\Big)ds \\
			-\bi \int_0^tS(t-s)&\Big(\widetilde 
			\Theta_R(u,k,s)\sum_{k\in \mathbb N^+}g(|u(s)|^2)|g'(|u(s)|^2) |u(s)|^2 u(s)
			Im(Q^{\frac 12}e_k)Q^{\frac 12}e_k \Big)ds,\\
			\Gamma_{Sto}^{R,k}  u(t):&=\mathbf i\int_0^t S(t-s) \Big(\widetilde \Theta_R(u,k,s)g(|u(s)|^2)u(s)\Big)dW^k(s),
			\end{align*}
			w}here $t\in [0,r]$, $W^k(s)=W(s+kr)-W(kr),$ $u\in \mathbb M_{\mathcal F_{kr}}^p(\Omega;C([0,r];\mathbb H^2))$
		and 
		$$\widetilde \Theta_R(u,k,s)=\theta_R(\|u^{R,k}\|_{C([0,kr];\mathbb H^2)}+\|u\|_{C([0,r];\mathbb H^2)}).$$
		
		For any different $v_1,v_2\in  \mathbb M_{\mathcal F_{kr}}^p(\Omega; C([0,r];\mathbb H^2))$, we define the stopping times $\tau_i=\inf\{t\in [0,T-kr]: \|u^{R,k}\|_{C([0,kr];\mathbb H^2)}+\|v_i\|_{C([0,t];\mathbb H^2)}\ge 2R\}\wedge r,$ $ i=1,2$ and assume that $\tau_1\le \tau_2$ for the convenience.  
		Then the same procedures yield that 
		this map is a strict contraction and has a fixed point $u^{R,k+1}$ for a small $r>0.$ Now, we define a process $u^R$ as 
		$u^{R}(t):=u^{R,k}(t)$ for $t\in [0,kr]$ and $u^{R}(t):=u^{R,k+1}(t)$ for $t\in [kr,kr+1].$ 
		It can be checked that $u^R$ satisfies \eqref{trun-reg} by the induction assumption and the definition of $\Gamma^R$. 
		Meanwhile, the uniqueness of the mild solution can be obtained by repeating the previous arguments.
	\hfill$\square$
	
	\begin{prop}\label{mul}
		Let the condition of Lemma \ref{trun} hold.
		There exists a unique local mild solution to \eqref{Reg-SlogS} with continuous $\mathbb H^2$-valued path. And the solution is defined 
		on a random interval $[0,\tau^*_{\epsilon}(u_0,\epsilon))$, where $\tau^*_{\epsilon}(u_0,\epsilon,\omega)$ is a stopping time such that $\tau^*_{\epsilon}(u_0,\epsilon,\omega)=+\infty$ or $\lim\limits_{t\to\tau^*_{\epsilon}}\|u^{\epsilon}(t)\|_{\mathbb H^2}=+\infty.$
	\end{prop}
	{\bf{Proof}}
		When $\widetilde g=1$, one can follow the same steps in the proof of \cite[Theorem 3.1]{BD03} to complete the proof. 
		 It suffices to consider the multiplicative noise case.
		Let $\{u^R\}_{R\in \mathbb N^+}\in \mathbb M_{\mathcal F}(\Omega; C([0,T];\mathbb H^2))$ be a sequence of solution constructed in Lemma \ref{trun}. And define a stopping time sequence $\tau_{R}:=\inf\{t\in[0,T]: \|u^R\|_{C([0,t];\mathbb H^2)}\ge R\}\wedge T.$
		Then $\tau_R>0$ is well-defined since $\|u^R\|_{C([0,t];\mathbb H^2)}$ is an increasing, continuous and $\mathcal F_t$-adapted process. 
		We claim that if $R_1\le R_2$, then $\tau_{R_1}\le \tau_{R_2}$ and $u^{R_1}=u^{R_2}, a.s.$ on $[0,\tau_{R_1}].$
		Let $\tau_{R_2,R_1}:=\inf\{t\in[0,T]: \|u^{R_2}\|_{C([0,t];\mathbb H^2)}\ge R_1\}\wedge T.$ 
		Then it holds that $\tau_{R_2,R_1}\le \tau_{R_2}$ and $\Theta_{R_2}(u^{R_2},t)=\Theta_{R_1}(u^{R_2},t)$ on $t\in [0,\tau_{R_2,R_1}].$ This implies that $\{(u^{R_2},\tau_{R_2,R_1})\}$ is a solution of \eqref{trun-reg} and that $u^{R_1}=u^{R_2}, a.s.$ on $\{t\le \tau_{R_2,R_1}\}.$ Thus we conclude that $\tau_{R_1}=\tau_{R_2,R_1}, a.s.$ and that $u^{R_1}=u^{R_2}$ for $\{t\le \tau_{R_1}\}$.
		
		Now consider the triple $\{u,(\tau_{R})_{R\in\mathbb N^+},\tau_{\infty}\}$ defined by $u(t):=u^R(t)$ for $t\in [0,\tau^R]$ and $\tau_{\infty}=\sup_{R\in\mathbb N^+}\tau_{R}.$ From Lemma \ref{trun}, we know that 
		$u\in \mathbb M_{\mathcal F}^p(\Omega, C([0,\tau];\mathbb H^2))$ satisfies \eqref{Reg-SlogS} for $t\le \tau_R$.
		The uniqueness of the local solution also holds.
		If we assume that $(u,\tau)$ and $(v,\sigma)$ are local mild solutions of \eqref{Reg-SlogS}, then $u(t)=v(t), a.s.$ on $\{t<\sigma \wedge \tau\}.$
		Let $R_1,R_2\in \mathbb N^+.$
		Set $\tau_{R_1,R_2}:=\inf\{t\in [0,T]: \max(\|u\|_{C([0,t];\mathbb H^2)}, \|v\|_{C([0,t];\mathbb H^2)})\ge R_1\}\wedge \sigma_{R_2} \wedge \sigma_{R_2}.$
		Then we have on $\{t\le\tau_{R_1,R_2}\},$ $(u,\tau_{R_1,R_2}),$ $(v,\tau_{R_1,R_2})$ are local mild solutions of \eqref{Reg-SlogS}. 
		The uniqueness 
		in Lemma \ref{trun} leads to $u=v$ on $\{t\le \tau_{R_1,R_2}\}.$ Letting $R_1,R_2\to \infty$, we complete the proof.\hfill$\square$

	If we assume that $g\in \mathcal C^{\mathbf s}_p(\mathbb R),$ $u_0\in \mathbb H^{\mathbf s}$, $\mathbf s\ge 2$, following the same procedures, we can also obtain the local existence of the solution $u_{\epsilon}$ in $\mathcal C([0,\tau^*_{\epsilon }]; \mathbb H^{\mathbf s})$ when $d\le 3$.
	When $d\ge 3$, one needs to use another regularization function $\log(\frac{|x|^2+\epsilon}{1+\epsilon |x|^2})$ and additional assumption on $g$. 
	In this case, we can get the local well-posedness in $\mathbb M_{\mathcal F}^p(\Omega; C([0,r];\mathbb H^1))$ based on Lemma \ref{sec-reg} in Appendix and previous arguments.  Since its proof is similar to that in Proposition \ref{mul}, we omit these details here and leave them to readers.  
	
	\begin{prop}\label{weak-loc}
		Let Assumption \ref{main-as} hold, and $f_{\epsilon}(|x|^2)$ $=\log(\frac{|x|^2+\epsilon}{1+\epsilon |x|^2})$. 
		Suppose that $d\in \mathbb N^+$, $u_0\in \mathbb H^2$ is $\mathcal F_0$-measurable and has any finite $p$th moment for $p\geq 1$. 
		Assume in addition that $Q$-Wiener process $W$ satisfies $\sum_{k\in \mathbb N^+}\|Q^{\frac 12}e_k\|_{L^2_\alpha}^2+\|Q^{\frac 12}e_k\|_{\mathbb H^1}^2<\infty$ when $\widetilde g=1$, and $\sum_{k\in \mathbb N^+}\|Q^{\frac 12}e_k\|_{W^{1,\infty}}^2+\|Q^{\frac 12}e_k\|_{\mathbb H^1}^2<\infty $ when $\widetilde g(x)=\mathbf i g(|x|^2)x$.
		Then there exists a unique local mild solution to \eqref{Reg-SlogS} with continuous $\mathbb H^1$-valued path. 
		 And the solution is defined 
			on a random interval $[0,\tau_{\epsilon}^*(u_0,\epsilon,\omega))$, where $\tau_{\epsilon}^*(u_0,\epsilon,\omega)$ is a stopping time such that $\tau^*(u_0,\epsilon,\omega)=+\infty$ or $\lim_{t\to\tau^*_{\epsilon}}\|u^{\epsilon}(t)\|_{\mathbb H^1}=+\infty.$
	\end{prop}

	\subsection{Global existence and uniform estimate of regularized SlogS equation}
	Due to the blow-up alternative results in section \ref{sec-loc}, it suffices to prove that $\sup_{t\in[0,\tau_{\epsilon}^*)}\|u^{\epsilon}(t)\|_{\mathbb H^2}<\infty,$ or $\sup_{t\in[0,\tau_{\epsilon}^*)}\|u^{\epsilon}(t)\|_{\mathbb H^1}<\infty, a.s.$ under corresponding assumptions. 
	In the following, we present several a priori estimates in strong sense to achieve our goal. 
	To simplify the presentation, we  omit some procedures 
	like mollifying the unbounded operator $\Delta$ and taking the limit on the regularization parameter. 
	More precisely, the mollifier $\widetilde \Theta_m$, $m\in \mathbb N^+$ may be defined by the Fourier transformation  (see e.g. \cite{BD03})
	\begin{align*}
	\mathbb F(\widetilde \Theta_m v)(\xi)=\widetilde \theta(\frac {|\xi|}{m})\widehat v(\xi),\; \xi\in \mathbb R^d,
	\end{align*}
	where $\widetilde \theta$ is a positive $\mathcal C^{\infty}$ function on $\mathbb R^+$, has a compact support satisfying 
	$\theta(x)=0,$ for $x\ge 2$ and $\theta(x)=1,$ for $0\le x\le 1.$
	Another choice of mollifier is via Yosida approximation $\Theta_{m}:=m(m-\Delta)^{-1}$ for $m\in \mathbb N^+$ (see e.g. \cite{Hor18}). 
	This kind of procedure is introduced to make that the It\^o formula can be applied rigorously to deducing several a priori estimates. 
	If $\mathcal O$ becomes a bounded domain equipped with periodic or homogenous Dirichlet boundary condition, the mollifier can be chosen as the Galerkin projection, and the approximated equation becomes the Galerkin approximation (see e.g. \cite{CH17,CHL16b,CHLZ17,CHS19}). 
	
	In this section, we assume that $u_0$ has finite $p$-moment for all $p\ge 1$ for simplicity. 
	We will also use Assumption \ref{main-as} to obtain the global existence of the mild solution of the regularized SlogS equation.
 When $d\le 3$ and $f_{\epsilon}(x)=\log(x+\epsilon),$ we assume that $\sum_{k}\|Q^{\frac 12}e_k\|_{\mathbb H^2}^2<\infty.$ When $d\in \mathbb N^+$ and $f_{\epsilon}(x)=\log(\frac {x+\epsilon}{1+\epsilon x}),$ we assume that $Q^{\frac 12}\in \mathbb L^1_2$ for the additive noise and $\sum_{k}\|Q^{\frac 12}e_k\|^2_{W^{1,\infty}}<\infty$ for the multiplicative noise.

	\subsection{A priori estimates in $L^2$ and $\mathbb H^1$}
	
	\begin{lm}\label{mass-pri}
		Let $T>0.$
		Under the condition of Proposition \ref{mul} or Proposition \ref{weak-loc}, assume that $(u^{\epsilon},\tau_{\epsilon}^*)$ be a local mild solution in $\mathbb H^1$. 
		Then for any $p\ge2,$ there exists a positive constant $C(Q,T,\lambda,p,u_0)>0$ such that
		\begin{align*}
		\E\Big[\sup_{t\in[0,\tau^*_{\epsilon}\wedge T)}\|u^{\epsilon}(t)\|^p\Big]\le 
		C(Q,T,\lambda,p,u_0).
		\end{align*}
	\end{lm}

	{\noindent\bf{Proof}}
		Take any stopping time $\tau<\tau^*_{\epsilon}\wedge T, a.s.$ Using the It\^{o} formula to 
		$M^k(u^{\epsilon}(t)),$ where $M(u^{\epsilon}(t)):=\|u^{\epsilon}(t)\|^2$ and $k\in \mathbb N^+$ or $k\ge 2$, we obtain that for $t\in [0,\tau]$ and the case $\widetilde g=1,$
		\begin{align*}
		&M^k(u^{\epsilon}(t))\\
		=&M^k(u_0^{\epsilon})
		+2k(k-1)\int_0^tM^{k-2}(u^{\epsilon}(s))\sum_{i\in \mathbb N^+} \<u^{\epsilon}(s),Q^{\frac 12}e_i\>^2ds
		\\
		&+k\int_{0}^tM^{k-1}(u^{\epsilon}(s))\sum_{i\in \mathbb N^+} \|Q^{\frac 12}e_i\|^2 ds
		+2k\int_{0}^tM^{k-1}(u^{\epsilon}(s))\<u^\epsilon(s),dW(s)\>,
		\end{align*} 
		and for $t\in [0,\tau]$ and the case $\widetilde g(x)=\bi g(|x|^2)x$, 
		\begin{align*}
		&\quad M^k(u^{\epsilon}(t))\\
		&=M^k(u_0^{\epsilon})+ {2k(k-1)}\int_0^tM^{k-2}(u^{\epsilon}(s))\sum_{i\in \mathbb N^+} \<u^{\epsilon}(s),\mathbf  i g(|u^{\epsilon}(s)|^2)u^{\epsilon}(s)Q^{\frac 12}e_i\>^2ds
		\\
		&\quad
		+2k\int_{0}^tM^{k-1}(u^{\epsilon}(s))\<u^\epsilon(s),\mathbf  i g(|u^{\epsilon}(s)|^2)u^{\epsilon}(s) dW(s)\>+2k\int_0^tM^{k-1}(u^{\epsilon}(s))\\
		&\quad\times\Big\<u^{\epsilon}(s), -\bi \sum_{k\in \N^+}Im(Q^{\frac 12}e_k) Q^{\frac 12}e_k g'(|u^{\epsilon}(s)|^2)g(|u^{\epsilon}(s)|^2)|u^{\epsilon}(s)|^2u^{\epsilon}(s)\Big\>ds.
		\end{align*}
		In particular, if  $W(t,x)$ is real-valued and $\widetilde g(x)=\bi g(|x|^2)x$, we have $M^k(u^{\epsilon}(t))=M^k(u^{\epsilon}_0),$ for $t\in[0,T]$ a.s.
		By  Assumption \ref{main-as} and conditions in Propositions \ref{mul} and \ref{weak-loc}, using the martingale inequality, the H\"older inequality and the Young inequality and Gronwall's inequality, we achieve that for all $k\ge1,$
		\begin{align*}
		\sup_{t\in [0,\tau]}\E\Big[M^k(u^{\epsilon}(t))\Big]\le C(T,k,u_0,{Q}).
		\end{align*}
		
		Next, taking the supreme over $t$ and repeating the above procedures,  we have that 
		\begin{align*}
		&\E\Big[\sup_{t\in [0,\tau]}M^k(u^{\epsilon}(t))\Big]\\
		\le& \E\Big[M^k(u_0^{\epsilon})\Big]+C(k)\E\Big[\int_0^\tau M^{k-2}(u^{\epsilon}(s))\sum_{i\in \mathbb N^+} (1+\|u^{\epsilon}(s)\|^2)\|Q^{\frac 12}e_i\|_{U}^2ds\Big]
		\\
		&\quad+C(k)\E\Big[\Big(\int_{0}^{\tau} M^{2k-2}(u^{\epsilon}(s))(1+\|u^{\epsilon}(s)\|^4) \sum_{i\in\N^+}\|Q^{\frac 12}e_i\|_{U}^2ds\Big)^{\frac 12}\Big],
		\end{align*}
		where $U=L^2$ for additive noise case and $U=L^{\infty}$ for multiplicative noise case.
		Applying the estimate of $\sup\limits_{t\in [0,\tau]}\E\Big[M^{k}(u^{\epsilon}(t))\Big]$, we complete the proof by taking $p=2k$. 
\hfill$\square$

	\begin{lm}\label{h1-pri}
		Let $T>0$. Under the condition of Proposition \ref{mul} or Proposition \ref{weak-loc}, assume that $(u^{\epsilon},\tau_{\epsilon}^*)$ is a local mild solution in $\mathbb H^1$.
		Then for any $p\ge2,$ there exists $C(Q,T,\lambda,p,u_0)>0$ such that
		\begin{align*}
		\E\Big[\sup_{t\in[0,\tau^*_{\epsilon}\wedge T)}\|u^{\epsilon}\|_{\mathbb H^1}^p\Big]\le 
		C(Q,T,\lambda,p,u_0).
		\end{align*}
	\end{lm}

{\noindent\bf{Proof}}	
		Take any stopping time $\tau<\tau^*_{\epsilon}\wedge T, a.s.$ 
		Applying the It\^{o} formula to 
		the kinetic energy  $K(u^{\epsilon}(t)):=\frac 12\|\nabla u^{\epsilon}(t)\|^2,$ for $k\in \mathbb N^+$ and using integration by parts,
		we obtain that for $t\in [0,\tau],$
		\begin{align*}
		&\quad K^k(u^{\epsilon}(t))\\
		&=K^k(u^{\epsilon}_0)
		+k\int_0^{t}K^{k-1}(u^{\epsilon}(s))\<\nabla u^{\epsilon}(s),\bi 2\lambda f'_{\epsilon}(|u^{\epsilon}(s)|^2)Re(\bar u^{\epsilon}(s)\nabla u^{\epsilon}(s)) u^{\epsilon}(s)\>ds\\
		&\quad+
		\frac 12k(k-1)\int_0^{t}K^{k-2}(u^{\epsilon}(s))\sum_{i\in \mathbb N^+}\<\nabla u^{\epsilon}(s), \nabla Q^{\frac 12}e_i\>^2ds\\
		&\quad 
		+\frac k 2\int_0^{t}K^{k-1}(u^{\epsilon}(s))\sum_{i\in \mathbb N^+}\<\nabla Q^{\frac 12}e_i,\nabla Q^{\frac 12}e_i\>ds\\
		&\quad 
		+k\int_0^{t}K^{k-1}(u^{\epsilon}(s))\<\nabla u^\epsilon(s),\nabla dW(s)\>
		\end{align*} 
		for the additive noise case,
		and {\small
			\begin{align*}
			&\quad K^k(u^{\epsilon}(t))\\
			&=K^k(u^{\epsilon}_0)+
			k\int_0^{t}K^{k-1}(u^{\epsilon}(s))\<\nabla u^{\epsilon}(s),\bi 2\lambda f'_{\epsilon}(|u^{\epsilon}(s)|^2)Re(\bar u^{\epsilon}(s)\nabla u^{\epsilon}(s)) u^{\epsilon}(s)\>ds\\
			&\quad+
			\frac 12 k(k-1)\int_0^{t}K^{k-2}(u^{\epsilon}(s))\sum_{i\in \mathbb N^+}\<\nabla u^{\epsilon}(s), \bi \nabla(g(|u^{\epsilon}(s)|^2)u^{\epsilon}(s) Q^{\frac 12}e_i)\>^2ds\\
			&\quad 
			+\frac k2\int_0^{t}K^{k-1}(u^{\epsilon}(s))\sum_{i\in \mathbb N^+}\<\nabla(g(|u^{\epsilon}(s)|^2)u^{\epsilon}(s) Q^{\frac 12}e_i),\nabla(g(|u^{\epsilon}(s)|^2)u^{\epsilon}(s) Q^{\frac 12}e_i)\>ds\\
			&\quad 
			+k \int_0^{t}K^{k-1}(u^{\epsilon}(s))\<\nabla u^{\epsilon}(s),\mathbf i \nabla g(|u^{\epsilon}(s)|^2)u^{\epsilon}(s) dW(s)\>\\
			&\quad+\frac k2 \int_0^{t}K^{k-1}(u^{\epsilon}(s))\<\Delta u^{\epsilon}(s),\sum_{i\in \N^+}|Q^{\frac 12}e_i|^2(g(|u^{\epsilon}(s)|^2))^2u^{\epsilon}(s)\>ds\\
			&\quad-k\int_0^{t}K^{k-1}(u^{\epsilon}(s))\<\Delta u^{\epsilon}(s),\bi \sum_{i\in \N^+}g(|u^{\epsilon}(s)|^2) g'(|u^{\epsilon}(s)|^2) |u^{\epsilon}(s)|^2 u^{\epsilon}(s)
			Im(Q^{\frac 12}e_i)Q^{\frac 12}e_i\>ds
			\end{align*} 
			f}or the multiplicative noise case. 
		Applying integration by parts, in the multiplicative noise case, we further obtain  
		\begin{align*}
		& K^k(u^{\epsilon}(t))\\
		\le& K^k(u^{\epsilon}_0)+
		k \int_0^{t}K^{k-1}(u^{\epsilon}(s))\<\nabla u^{\epsilon},\bi2 \lambda f'_{\epsilon}(|u^{\epsilon}(s)|^2)Re(\bar u^{\epsilon} \nabla u^{\epsilon})u^{\epsilon}(s)\>ds\\
		&+
		C_k \int_0^{t}K^{k-2}(u^{\epsilon}(s))\sum_{i\in \mathbb N^+}\Big(\<\nabla u^{\epsilon},  g(|u^{\epsilon}|^2)u^{\epsilon} \nabla Q^{\frac 12}e_i)\>^2+\<\nabla u^{\epsilon}, g'(|u^{\epsilon}|^2)Re(\bar u^{\epsilon}\nabla u^{\epsilon})u^{\epsilon} Q^{\frac 12}e_i)\>^2\Big)ds\\
		& 
		+C_k\int_0^{t}K^{k-1}(u^{\epsilon}(s))\sum_{i\in \mathbb N^+}\Big(\|g(|u^{\epsilon}|^2) u^{\epsilon} \nabla Q^{\frac 12}e_i\|^2+
	\\
		&\quad+\|g(|u^{\epsilon}|^2)\nabla u^{\epsilon} Q^{\frac 12}e_i\|^2+\|g'(|u^{\epsilon}|^2) Re(\bar u^{\epsilon}\nabla u^{\epsilon})u^{\epsilon} Q^{\frac 12}e_i\|^2\Big)ds\\
		&
		+k \int_0^{t}K^{k-1}(u^{\epsilon}(s))\<\nabla u^\epsilon,\mathbf i \nabla g(|u^{\epsilon}|^2)u^{\epsilon} dW(s)\>\\
		&+C_k\int_0^{t}K^{k-1}(u^{\epsilon}(s))\sum_{i\in \mathbb N^+}\Big(|\<\nabla u^{\epsilon},(\nabla Im(Q^{\frac 12}e_i) Q^{\frac 12}e_i+\nabla Q^{\frac 12}e_i Im(Q^{\frac 12}e_i)) g'(|u^{\epsilon}|^2) g(|u^{\epsilon}|^2) |u^{\epsilon}|^2u^{\epsilon}\>|\\
		&\quad+|\<\nabla u^{\epsilon}, Im(Q^{\frac 12}e_i)Q^{\frac 12}e_i g''(|u^{\epsilon}|^2) g(|u^{\epsilon}|^2) |u^{\epsilon}|^2u^{\epsilon} Re(\bar u^{\epsilon}\nabla u^{\epsilon})\>| \\
		&\quad +|\<\nabla u^{\epsilon}, Im(Q^{\frac 12}e_i)Q^{\frac 12}e_i|g'(|u^{\epsilon}|^2)|^2 |u^{\epsilon}|^2u^{\epsilon} Re(\bar u^{\epsilon}\nabla u^{\epsilon})\>| \\
		&\quad +|\<\nabla u^{\epsilon}, Im(Q^{\frac 12}e_i)Q^{\frac 12}e_ig'(|u^{\epsilon}|^2) g(|u^{\epsilon}|^2)  (2Re(\bar u^{\epsilon}\nabla u^{\epsilon})u^{\epsilon}+|u^{\epsilon}|^2\nabla u^{\epsilon})\>|\Big)ds.
		\end{align*} 
		By using the property of $g$ in Assumption \ref{main-as} and conditions on $Q$, and applying H\"older's, Young's and Burkholder's  inequalities, we achieve that for small $\epsilon_1>0,$
		\begin{align*}
		&\quad\E[\sup_{r\in[0,t]}K^k(u^{\epsilon}(r))]\\
		&\le \E[K^k(u^{\epsilon}_0)]+
		C_k |\lambda| \E\Big[ \int_{0}^{t}\sup_{r\in[0,s]} K^{k-1}(u^{\epsilon}(r)) \|\nabla u^{\epsilon}(r)\|^2 ds\Big]
		\\
		&\quad+C_kT\E\Big[\sup_{r\in[0,t]}K^{k-2}(u^{\epsilon}(r))\sum_{i\in \mathbb N^+}\|\nabla u^{\epsilon}(r)\|^2(1+\|\nabla u^{\epsilon}(r)\|^2+\|u^{\epsilon}(r)\|^2) \|Q^{\frac 12}e_i\|_{U}^2\Big]\\
		&\quad+ C_k\E\Big[\Big(\int_0^{t}K^{2k-2}(u^{\epsilon}(s))\sum_{i\in\mathbb N^+}\|\nabla u^{\epsilon}(r)\|^2(1+\|\nabla u^{\epsilon}(r)\|^2+\|u^{\epsilon}(r)\|^2 \|\nabla Q^{\frac 12}e_i\|_U^2ds\Big)^{\frac 12}\Big]\\
		&\le \E[K^k(u^{\epsilon}_0)]+C(T,\lambda,k,\epsilon_1,Q)
		+\epsilon_1 \E[\sup_{r\in[0,t]}K^k(u^{\epsilon}(r))]\\
		&\quad+C_k \E\Big[ \int_{0}^{t}\sup_{r\in[0,s]} K^{k}(u^{\epsilon}(r)) dr\Big],
		\end{align*} 
		where $\|Q^{\frac 12}e_i\|_U^2=\|Q^{\frac 12}e_i\|^2_{\mathbb H^1}$ for the additive noise case and $\|Q^{\frac 12}e_i\|_U^2=\|Q^{\frac 12}e_i\|^2_{\mathbb H^1}+\|Q^{\frac 12}e_i\|^2_{W^{1,\infty}}$ for the multiplicative noise case. 
		Applying Gronwall's inequality, we complete the proof by taking $p=2k$. 
	\qed
	
	From the above proofs of Lemmas \ref{mass-pri} and \ref{h1-pri},
	it is not hard to see that to obtain $\epsilon$-independent estimates, the boundedness restriction $\sup_{x\ge 0}|g(x)|<\infty$ may be not necessary in the case that $W(t,x)$ is real-valued. 
		We present such result in the following which is the key of the global well-posedness of an SlogS equation with super-linear growth diffusion in next section. 
	
	\begin{lm}\label{lm-con-h1}
		Let $T>0$ and $(u^{\epsilon},\tau^*_{\epsilon})$ be a local mild solution in $\mathbb H^{\bs}$, $\bs\ge 1$ for any $p\ge 1.$ Assume that $u_0\in \mathbb H^1\cap L^2_{\alpha}$, for some $\alpha\in (0,1]$, is $\mathcal F_0$-measurable and has any finite $p$th moment, and $W(t,x)$ is real-valued with $\sum_{i} \|Q^{\frac 12}e_i\|_{\mathbb H^1}^2+\|Q^{\frac 12}e_i\|_{W^{1,\infty}}^2<\infty.$  Let $\widetilde g(x)=\mathbf ig(|x|^2)x,$ $g \in \mathcal C^1_b(\mathbb R)\cap \mathcal C(\mathbb R)$ satisfy the growth condition and  the embedding condition,
		\begin{align*}
		&\sup_{x\in [0,\infty)}|g'(x)x|\le C_g,\\
		&\|vg(|v|^2)\|_{L^q}\le C_d (1+\|v\|_{\mathbb H^1}+\|v\|_{L^2_{\alpha}}),
		\end{align*}
		for some $q\ge 2$, where $C_g>0$ depends on $g$, and $C_d>0$ depends on $\mathcal O$, $d$, $\|v\|$. 
		Then it holds that $M(u^{\epsilon}_t)=M(u_0)$ for $t\in [0,\tau_{\epsilon}).$
		Furthermore, there exists a positive constant $C(Q,T,\lambda,p,u_0)$ such that 
		\begin{align*}
		\E\Big[\sup_{t\in[0,\tau^*_{\epsilon}\wedge T)}\|u^{\epsilon}(t)\|_{\mathbb H^1}^p\Big]\le 
		C(Q,T,\lambda,p,u_0).
		\end{align*}
	\end{lm}
	\textbf{Proof}
		The proof is similar to those of Lemmas \ref{mass-pri} and \ref{h1-pri}. We only need to modify the estimation involved with $g.$
		The mass conservation is not hard to be obtained since the calculations in Lemmas \ref{mass-pri} only use the assumptions that $g\in \mathcal C(\mathbb R)$ and $W(t,x)$ is real-valued. Therefore, we focus on the estimate in $\mathbb H^1.$ 
		We only show estimation about $\E\Big[K^k(u^{\epsilon}(t))\Big]$ since the proof on $\E\Big[\sup_{t\in[0,\tau^*_{\epsilon}\wedge T)}K^k(u^{\epsilon}(t))\Big]$ is similar.
		Then following the same steps in Lemma \ref{h1-pri}, we get that for $\frac 1q+\frac 1{q'}=\frac 12$, 
		\begin{align*}
		&\E\Big[K^k(u^{\epsilon}(t))\Big]\\
		\le& \E\Big[K^k(u^{\epsilon}_0)\Big]+
		k \E\Big[\int_0^{t}K^{k-1}(u^{\epsilon}(s))\<\nabla u^{\epsilon},\bi2 \lambda f'_{\epsilon}(|u^{\epsilon}(s)|^2)Re(\bar u^{\epsilon} \nabla u^{\epsilon})u^{\epsilon}(s)\>ds\Big]\\
		&+
		\E\Big[C_k\int_0^{t}K^{k-2}(u^{\epsilon}(s))\sum_{i\in \mathbb N^+}\Big(\<\nabla u^{\epsilon},  g(|u^{\epsilon}|^2)u^{\epsilon} \nabla Q^{\frac 12}e_i)\>^2+\<\nabla u^{\epsilon}, g'(|u^{\epsilon}|^2)Re(\bar u^{\epsilon}\nabla u^{\epsilon})u^{\epsilon} Q^{\frac 12}e_i)\>^2\Big)ds\Big]\\
		& 
		+C_k\E\Big[\int_0^{t}K^{k-1}(u^{\epsilon}(s))\sum_{i\in \mathbb N^+}\Big(\|g(|u^{\epsilon}|^2) u^{\epsilon} \nabla Q^{\frac 12}e_i\|^2+\|g'(|u^{\epsilon}|^2) Re(\bar u^{\epsilon}\nabla u^{\epsilon})u^{\epsilon} Q^{\frac 12}e_i\|^2\Big)ds\Big]\\
		\le& \E\Big[K^k(u^{\epsilon}_0)\Big]+
		C_k \E\Big[\int_0^{t}K^{k}(u^{\epsilon}(s))ds\Big]\\
		&+
		\E\Big[C_k\int_0^{t}K^{k-1}(u^{\epsilon}(s))\sum_{i\in \mathbb N^+}\Big(\|g(|u^{\epsilon}|^2)u^{\epsilon}\|_{L^q}^2 \|\nabla Q^{\frac 12}e_i\|_{L^{q'}}^2+ \|\nabla u^{\epsilon}\|^2\|Q^{\frac 12}e_i\|^2\Big)ds\Big].
		\end{align*} 
		Applying  the embedding condition on $g$, mass conservation law and the procedures in the proof of Proposition \ref{wei-hol-1}, we complete the proof by taking $p=2k$ and Gronwall's inequality.
	\qed
	
	The embedding condition is depending on the assumption on $Q$ and $d$. One example which satisfies the embedding condition on $g$ and is not bounded is $g(x)=\log(c+x),$ for $x\ge 0$ and $c>0.$ Let us verify this example on $\mathcal O=\mathbb R^d.$ 
	If the domain $\mathcal O$  is  bounded, one can obtain the similar estimate. 
	We apply the Gagliardo--Nirenberg interpolation inequality and get that for $q>2$ and small enough $\eta>0,$
	\begin{align*}
	\|g(|v|^2)v\|_{L^q}^q&\le  \int_{c+|v|^2\ge 1}|\log(c+|v|^2)|^q|v|^qdx+ \int_{c+|v|^2\le 1}|\log(c+|v|^2)|^q|v|^qdx\\
	&\le C\Big(\|v\|_{L^q}^q+\|v\|_{L^{q+q\eta}}^{q+q\eta}+\|v\|_{L^{q-q\eta}}^{q-q\eta}\Big)\\
	&\le  C\Big(\|v\|^{1-\alpha_0}\|\nabla v\|^{\alpha_0}+(\|v\|^{1-\alpha_1}\|\nabla v\|^{\alpha_1})^{1+\eta}+(\|v\|^{1-\alpha_2}\|\nabla v\|^{\alpha_2})^{1-\eta}\Big)^q\\
	&\le C \Big({\|\nabla v\|^{q\alpha_0}}+\|\nabla v\|^{(q+q\eta)\alpha_1}+\|\nabla v\|^{(q-q\eta)\alpha_2}\Big),
	\end{align*}
	where $\alpha_0=\frac {d(q-2)}{2q}, \alpha_1=\frac {d(q(1+\eta)-2)}{2q(1+\eta)}$ and $\alpha_2=\frac {d(q(1-\eta)-2)}{2q(1-\eta)}$ satisfies $\alpha_i\in (0,1), i=0,1,2.$
	When $q=2,$ similar calculations, together with the interpolation inequality in Lemma \ref{wei-hol},  yield that
	\begin{align*}
	\|g(|v|^2)v\|^2
	&\le C\Big(\|v\|^2+\|v\|_{L^{2+2\eta}}^{2+2\eta}+\|v\|_{L^{2-2\eta}}^{2-2\eta}\Big)\\
	&\le C \Big(\|\nabla v\|^{2}+\|\nabla v\|^{(2+2\eta)\alpha_1}+\|v\|_{L^2_{\alpha_2}}^{\frac {d\eta}{\alpha_2}}\Big),
	\end{align*}
 where $\alpha_1=\frac {d(2(1+\eta)-2)}{4(1+\eta)}\in(0,1)$ and $\alpha_2\in(\frac {d\eta}{2-2\eta},1).$

	\subsection{$L^2_{\alpha}$-estimate and modified energy}
	
	Beyond the $L^2$ and $\mathbb H^1$ estimates, we also need the uniform boundedness in $L_{\alpha}^2,$ $\alpha>0$ to show the strong convergence of  $\{u^{\epsilon}\}_{\epsilon>0}$ when $\mathcal O=\mathbb R^d.$ 
	We would like to mention that when $\mathcal O$ is a bounded domain, such estimate in $L_{\alpha}^2,$ $\alpha>0$ is not necessary. To this end, the following useful weighted interpolation inequality is introduced. 
	
	\begin{lm}\label{wei-hol}
		Let $d\in \mathbb N^+$ and $\eta\in (0,1)$. 
		Then for $\alpha>\frac {d\eta}{2-2\eta}$,  it holds that for some $C=C(d)>0$,
		\begin{align*}
		\|v\|_{L^{2-2\eta}}\le C\|v\|^{1-\frac {d\eta}{2\alpha(1-\eta)}}\|v\|_{L^2_{\alpha}}^{\frac {d\eta}{2\alpha(1-\eta)}},\; v\in L^2\cap L^2_{\alpha}.
		\end{align*} 
	\end{lm}
	\textbf{Proof}
		Using the Cauchy-Schwarz inequality and $\alpha>\frac {d\eta}{2-2\eta}$, we have that for any $r>0,$
		\begin{align*}
		\|v\|_{L^{2-2\eta}}^{2-2\eta}&\le \int_{|x|\le r}|v(x)|^{2-2\eta}dx+\int_{|x|\ge r}\frac {|x|^{ \alpha({2-2\eta})}|v(x)|^{2-2\eta}}{|x|^{\alpha({2-2\eta})}}dx\\
		&\le Cr^{ d \eta} \|v\|^{ {2-2\eta}}+C\|v\|_{L^2_{\alpha}}^{{2-2\eta}}(\int_{|x|\ge r}\frac 1{|x|^{\frac {\alpha(2-2\eta)}{\eta}}}dx)^{\eta}\\
		&\le Cr^{ d \eta} \|v\|^{ {2-2\eta}}+C r^{- \alpha(2-2\eta)+d\eta} \|v\|_{L^2_{\alpha}}^{{2-2\eta}}.
		\end{align*}
		Let $r=(\frac {\|v\|_{L^2_{\alpha}}}{\|v\|})^{\frac {1} {\alpha}}$, we complete the proof.
	\qed

	\begin{prop}\label{wei-hol-1}
		Let $T>0$, $\mathcal O=\mathbb R^d$, $d\in \mathbb N^+$.
		Under the condition of Proposition \ref{mul} or Proposition \ref{weak-loc}, let $(u^{\epsilon},\tau_{\epsilon}^*)$ be a local mild solution in $\mathbb H^{\bs}$, $\bs\ge 1$ for any $p\ge 1.$ 
		Let $u_0\in  L^2_{\alpha}\cap \mathbb H^1$, for some $\alpha\in (0,1]$, have any finite $p$th moment.
		Then the solution $u^{\epsilon}$ of regularized problem satisfies for $\alpha\in (0,1],$
		\begin{align*}
		\E\Big[\sup_{t\in[0,\tau_{\epsilon}^*\wedge T)}\|u^{\epsilon}(t)\|_{L^2_{\alpha}}^{2p}\Big]&\le C(Q,T,\lambda,p,u_0).
		\end{align*}
	\end{prop}
	
	{\noindent\bf{Proof}} 
		We first introduce the stopping time $$\tau_R=\inf\Big\{t\in [0,T]: \sup_{s\in[0,t]}\|u^{\epsilon}(s)\|_{L^2_{\alpha}}\ge R\Big\}\wedge \tau_{\epsilon}^*,$$ then  show that $\E\Big[\sup\limits_{t\in[0,\tau_{R}]}\|u^{\epsilon}(t)\|_{L^2_{\alpha}}^{2p}\Big]\le C(T,u_0,Q,p)$ independent of $R$. After taking $R\to \infty$, we get $\tau_R=\tau_{\epsilon}^*, a.s.$ For simplicity, we only prove uniform  upper bound when {\color{red}$p=1$.}
		
		Taking $0<t\le t_1\le \tau_R,$
		and applying the It\^o formula to $\|u^{\epsilon}\|_{L^2_{\alpha}}^{2}=\int_{\mathbb R^d}(1+|x|^2)^{\alpha}|u^{\epsilon}|^2dx$,
		we get 
		\begin{align*}
		\|u^{\epsilon}(t)\|_{L^2_{\alpha}}^{2}
		&=\|u^{\epsilon}_0\|_{L^2_{\alpha}}^{2}
		+\int_0^t 2 \<(1+|x|^2)^{\alpha}u^{\epsilon}(s),\bi \Delta u^{\epsilon}(s)\>ds+\int_0^t 2 \sum_{i\in \mathbb N^+}\<(1+|x|^2)^{\alpha}Q^{\frac 12}e_i, Q^{\frac 12}e_i\>ds\\
		&\quad+\int_0^t 2\<(1+|x|^2)^{\alpha}u^{\epsilon}(s), \bi f_{\epsilon}(|u^{\epsilon}(s)|^2)u^{\epsilon}(s)\>ds
		+\int_0^t 2\<(1+|x|^2)^{\alpha}u^{\epsilon}(s),dW(s)\>
		\end{align*}
		for additive noise case,
		and 
		\begin{align*}
		\|u^{\epsilon}(t)\|_{L^2_{\alpha}}^{2}
		&=\|u^{\epsilon}_0\|_{L^2_{\alpha}}^{2}
		+2\int_0^t \<(1+|x|^2)^{\alpha}u^{\epsilon}(s),\bi \Delta u^{\epsilon}(s)+\bi f_{\epsilon}(|u^{\epsilon}(s)|^2)u^{\epsilon}(s)\>ds\\
		&\quad-\int_0^t \<(1+|x|^2)^{\alpha}u^{\epsilon}(s)(s),  (g(|u^{\epsilon}(s)|^2))^2u^{\epsilon}(s)\sum_{i}|Q^{\frac 12}e_i|^2\>ds\\
		&\quad+\int_0^t \sum_{i}\<(1+|x|^2)^{\alpha} g(|u^{\epsilon}(s)|^2)u^{\epsilon}(s)Q^{\frac 12}e_i,  g(|u^{\epsilon}(s)|^2)u^{\epsilon}(s)Q^{\frac 12}e_i\>ds\\
		&\quad -2\int_0^t \<(1+|x|^2)^{\alpha}u^{\epsilon}(s),  \bi \sum_{k\in \N^+}g(|u^{\epsilon}(s)|^2) g'(|u^{\epsilon}(s)|^2) |u^{\epsilon}(s)|^2 u^{\epsilon}(s)
		Im(Q^{\frac 12}e_k)Q^{\frac 12}e_k \>ds\\
		&\quad +2\int_0^t\<(1+|x|^2)^{\alpha}u^{\epsilon}(s),\mathbf i g(|u^{\epsilon}(s)|^2)u^{\epsilon}(s)dW(s)\>
		\end{align*}
		for multiplicative noise case. 
Using integration by parts, then taking supreme over $t\in [0,t_1]$ and applying the Burkerholder inequality, we deduce 
		\begin{align*}
		\E\Big[\sup_{t\in [0,t_1]}\|u^{\epsilon}(t)\|_{L^2_{\alpha}}^{2}\Big]
		&\le \E\Big[\|u^{\epsilon}_0\|_{L^2_{\alpha}}^{2}\Big]
		+C_{\alpha} \E \Big[\int_0^T\Big| \<(1+|x|^2)^{\alpha-1} x u^{\epsilon}(s),\nabla u^{\epsilon}(s)\> \Big|ds\Big]\\
		&\quad+C\Big(\E \Big[\int_0^{t_1}\sum_{i\in\mathbb N^+}\|(1+|x|^2)^{\frac \alpha 2}{{u^{\epsilon}(s)}}\|^2 \|(1+|x|^2)^{\frac \alpha 2} Q^{\frac 12}e_i\|^2 ds\Big]\Big)^{\frac 12}
		\end{align*}
		for additive noise case,
		and 
		\begin{align*}
		\E\Big[\sup_{t\in [0,t_1]}\|u^{\epsilon}(t)\|_{L^2_{\alpha}}^{2}\Big]
		&\le \E\Big[\|u^{\epsilon}_0\|_{L^2_{\alpha}}^{2}\Big]
		+C_{\alpha} \E \Big[\int_0^T\Big| \<(1+|x|^2)^{\alpha-1} x u^{\epsilon}(s),\nabla u^{\epsilon}(s)\> \Big|ds\Big]\\
		&\quad+C\Big(\E \Big[\int_0^{t_1}\sum_{i\in\mathbb N^+}\|(1+|x|^2)^{\frac \alpha 2}u^{\epsilon}(s)\|^4 \|Q^{\frac 12}e_i\|_{L^{\infty}}^2 ds\Big]\Big)^{\frac 12}
		\end{align*}
		for multiplicative noise case.
		Then Young's and Gronwall's  inequalities, together with a priori estimate of $u^{\epsilon}$ in $\HH^1$, lead to the desired result.
	\qed
	
	\begin{cor}\label{cor-con-l2a}
		Under the condition of Lemma \ref{lm-con-h1}, the solution $u^{\epsilon}$ of regularized problem satisfies for $\alpha\in (0,1],$
		\begin{align*}
		\E\Big[\sup_{t\in[0,\tau_{\epsilon}^*\wedge T)}\|u^{\epsilon}(t)\|_{L^2_{\alpha}}^{2p}\Big]&\le C(Q,T,\lambda,p,u_0).
		\end{align*}
	\end{cor}

	It is not possible to obtain the uniform bound of the exact solution in $L^2_{\alpha}$ for $\alpha\in (1,2]$ like  the deterministic case. The main reason is that the rough driving noise leads to low H\"older regularity in time and loss of uniform estimate in $\mathbb H^2$ for the mild solution. We can not expect that the mild solution of \eqref{Reg-SlogS} enjoys  $\epsilon$-independent estimate in $\mathbb H^2.$ More precisely, we prove that  applying the regularization $f_{\epsilon}(|x|^2)=\log(|x|^2+\epsilon)$ in Proposition \ref{mul}, one can only expect $\epsilon$-dependent estimate in $\mathbb H^2$.  We omit the tedious calculation and procedures, and present a sketch of the proof for Lemma \ref{h2-pri} and Propostion \ref{wei-hol-2} in Appendix.
	
	\begin{lm}\label{h2-pri}
		Let $T>0$ and $d=1$. Under the condition of Proposition \ref{mul}, assume that $(u^{\epsilon},\tau_{\epsilon}^*)$ is the local mild solution in $\mathbb H^1$. In addition assume that $\sum_{i\in \mathbb N^+}\|Q^{\frac 12}e_i\|_{W^{2,\infty}}^2<\infty$ for the multiplicative noise. Then for any $p\ge2,$ there exists a positive $C(Q,T,\lambda,p,u_0)$ such that
		\begin{align*}
		\E\Big[\sup_{t\in[0,\tau_{\epsilon}^*\wedge T)}\|u^{\epsilon}(t)\|_{\mathbb H^2}^{2p}\Big]\le 
		C(Q,T,\lambda,p,u_0)(1+\epsilon^{-2p}).
		\end{align*}
	\end{lm}
	
	\begin{prop}\label{wei-hol-2}
		Assume that $\mathcal O=\mathbb R^d$. Let $T>0$ and $d=1$. Under the condition of Proposition \ref{mul}, assume that $(u^{\epsilon},\tau_{\epsilon}^*)$ is the local mild solution in $\mathbb H^1$. 
		Let $u_0\in  L^2_{\alpha}$, for some $\alpha\in (1,2]$. In addition assume that  $\sum_{i}\|Q^{\frac 12}e_i\|_{L^2_{\alpha}}^2<\infty$ in the additive noise case and that $\sum_{i\in \mathbb N^+}\|Q^{\frac 12}e_i\|_{W^{2,\infty}}^2<\infty$ for the multiplicative noise. Then the solution $u^{\epsilon}(t)$ of regularized problem satisfies for $\alpha\in (1,2],$
		\begin{align*}
		\E\Big[\sup_{t\in[0,\tau_{\epsilon}^*\wedge T)}\|u^{\epsilon}(t)\|_{L^2_{\alpha}}^{2p}\Big]&\le C(Q,T,\lambda,p,u_0)(1+\epsilon^{-2p}).
		\end{align*}
	\end{prop}
	
 The above results indicate that both spatial and temporal regularity for SLogS equation are rougher than deterministic LogS equation. 
	
	In the following, we present the behavior of the regularized energy for the RSlogS equation.
	When applying $f_{\epsilon}(|x|^2)=\log(|x|^2+\epsilon),$  the modified energy of \eqref{Reg-SlogS} becomes $$H
	_{\epsilon}(u^{\epsilon}(t)):=K(u^{\epsilon})-\frac \lambda{2}F_{\epsilon}(|u^{\epsilon}|^2)$$ with $F_{\epsilon}(|u^{\epsilon}|^2)=\int_{\mathcal O}\Big((\epsilon+|u^{\epsilon}|^2)\log(\epsilon+|u^{\epsilon}|^2)-|u^{\epsilon}|^2-\epsilon\log\epsilon\Big) dx.$
	When using another regularization function $f_{\epsilon}(|u|^2)=\log(\frac {\epsilon+|u|^2}{1+\epsilon|u|^2}),$ its regularized entropy in the modified energy $H_{\epsilon}$ becomes $$ F_{\epsilon}(|u^{\epsilon}|^2)=\int_{\mathcal O}\Big(|u^{\epsilon}|^2\log(\frac {|u^{\epsilon}|^2+\epsilon}{1+|u^{\epsilon}|^2\epsilon})+\epsilon \log(|u^{\epsilon}|^2+\epsilon)-\frac1{\epsilon}\log(\epsilon |u^{\epsilon}|^2+1)-\epsilon\log(\epsilon)\Big) dx.$$
In general, the modified energy is defined by the regularized entropy $\widetilde F(\rho)=\int_\mathcal O\int_0^{\rho}f_{\epsilon}(s)dsdx,$ where  $f_{\epsilon}(\cdot)$ is a suitable approximation of $\log(\cdot).$
	We remark that the regularized energy is well-defined when $\mathcal O$ is a bounded domain. The additional constant term $\epsilon\log(\epsilon)$ ensures that the regularized energy is still well-defined when $\mathcal O=\mathbb R^d.$  We leave the proof of Proposition \ref{eng-pri} in Appendix.

	\begin{prop}\label{eng-pri}
		Let $T>0$. Under the condition of Proposition \ref{mul} or Proposition \ref{weak-loc}, assume that $(u^{\epsilon},\tau_{\epsilon}^*)$ is the local mild solution in $\mathbb H^1$.
		Then for any $p\ge2,$ there exists a positive constant $C(Q,T,\lambda,p,u_0)$ such that
		\begin{align*}
		\E\Big[\sup_{t\in[0,\tau_{\epsilon}\wedge T)}|H_{\epsilon}(u^{\epsilon}(t))|^p\Big]\le 
		C(Q,T,\lambda,p,u_0).
		\end{align*}
	\end{prop}
	
	Below, we present the global existence of the unique mild solution for Eq. \eqref{Reg-SlogS} in $\mathbb H^2$ based on Proposition \ref{mul}, Lemma \ref{h2-pri}, Proposition \ref{weak-loc} and  Lemma \ref{h1-pri}, as well as a standard argument in \cite{BD03}.
	
	\begin{prop}\label{prop-well-reg}
		Let Assumption \ref{main-as} hold and $(u^{\epsilon},\tau_{\epsilon}^*)$ be the local mild solution in Proposition \ref{mul} or Proposition \ref{weak-loc}.
		Then the mild solution $u^{\epsilon}$ in $\mathbb H^1$ is global, i.e., $\tau^{*}_{\epsilon}=+\infty, a.s.$ 
		In addition assume that $d=1$ in Proposition \ref{mul}, the mild solution $u^{\epsilon}$ in $\mathbb H^2$ is global.
	\end{prop}
	
	\section{Well-posedness for SLogS equation}
	Based on the a priori estimates of the regularized problem, we are going to prove the strong convergence of any sequence of the solutions of the regularized problem.
	This immediately implies that the existence and uniqueness of the mild solution in $L^2$
	for SLogS equation. 
	\subsection{Well-posedness for SLogS equation via strong convergence approximation}
	In this part, we not only show the strong convergence of a sequence of  solutidregularized problems, but also give the explicit strong convergence rate. 
	The strong convergence rate of the regularized SLogS equation will make a great contribution to the numerical analysis of numerical schemes for the SLogS equation. And this topic will be studied in a companion paper.
	For the strong convergence result, we only present the mean square convergence rate result since the proof of the strong convergence rate in $L^q(\Omega), q\ge 2$ is similar. 
	In this section, the properties of regularization function $f_{\epsilon}$ in Lemmas \ref{frist-loc} and \ref{sec-reg} will be frequently used.

	In the multiplicative noise case,  Assumption \ref{ass-strong} is needed to obtain the strong convergence rate of the solution of  Eq. \eqref{Reg-SlogS}. 
	We remark that the assumption can be weaken if one only wants to obtain the strong convergence instead of deriving a convergence rate.
	Some sufficient condition for \eqref{con-g1} in Assumption \ref{ass-strong} is 
	$$\Big|(g'(|x|^2)g(|x|^2)|x|^2-g'(|y|^2)g(|y|^2)|y|^2)(|x|^2-|y|^2)\Big|\le C_g|x-y|^2, x,y\in  \mathbb C$$
	or 
	$$
	\Big|(g'(|x|^2)g(|x|^2)|x|^2-g'(|y|^2)g(|y|^2)|y|^2)(|x|+|y|)\Big|\le C_g|x-y|, x,y\in  \mathbb C.
	$$
	Functions like $1, \frac 1{c+x},\frac x{c+x},  \frac x{c+x^2}, \log(\frac {c+|x|^2}{1+c|x|^2})$  with $c>0$, etc., will satisfy Assumption \ref{ass-strong}.

	The main idea of the proof lies on showing that for a decreasing sequence $\{\epsilon_{n}\}_{n\in \mathbb N^+}$
	satisfying $\lim\limits_{n\to\infty}\epsilon_n=0$, $\{u^{\epsilon_n}\}_{n\in \mathbb N^+}$ must be a Cauchy sequence in $C([0,T];L^p(\Omega; \mathbb H)),$ $p\ge 2.$ As a consequence, we obtain that there exists a limit process $u$ in $C([0,T];L^p(\Omega; \mathbb H))$ which is shown to be independent of the sequence $\{u^{\epsilon_n}\}_{n\in \mathbb N^+}$ and is the unique 
mild solution of the mild form of \eqref{SlogS}.

	\textbf{[Proof of Theorem \ref{mild-general}]}
		Based on Proposition \ref{prop-well-reg}, we can construct a sequence of mild solutions $\{u^{\epsilon_n}\}_{n\in \mathbb N^+}$ of Eq. \eqref{Reg-SlogS} with $f_{\epsilon^n}(|x|^2)=\log(\frac {\epsilon_n+|x|^2}{1+\epsilon_n |x|^2})$. Here the decreasing sequence $\{\epsilon_{n}\}_{n\in \mathbb N^+}$ satisfies $\lim\limits_{n\to\infty}\epsilon_n=0.$ We use the following steps to complete the proof. For simplicity, we only present the details for $p=2$ since the procedures for $p>2$ are similar.
		
		Step 1: $\{u^{\epsilon_n}\}_{n\in \mathbb N^+}$ is a Cauchy sequence in $L^2(\Omega; C([0,T];\mathbb H)).$
		
		Fix different $n,m\in \mathbb N^+$ such that $n<m.$ 
		Subtracting the equation of $u^{\epsilon_n}$ from the equation of $u^{\epsilon_m}$, we have that 
		\begin{align*}
		d(u^{\epsilon_m}-u^{\epsilon_n})=\bi \Delta (u^{\epsilon_m}-u^{\epsilon_n})dt
		+\bi \lambda (f_{\epsilon_m}(|u^{\epsilon_m}|^2)u^{\epsilon_m}-f_{\epsilon_n}(|u^{\epsilon_n}|^2)u^{\epsilon_n})dt
		\end{align*}
		for additive noise case,
		and 
		\begin{align*}
		&\quad d(u^{\epsilon_m}-u^{\epsilon_n})\\
		&=\bi \Delta (u^{\epsilon_m}-u^{\epsilon_n})dt	 
		+\bi \lambda (f_{\epsilon_m}(|u^{\epsilon_m}|^2)u^{\epsilon_m}-f_{\epsilon_n}(|u^{\epsilon_n}|^2)u^{\epsilon_n})dt\\
		&\quad -\frac 12\sum_{k\in\mathbb N^+}|Q^{\frac 12}e_k|^2\Big(|g(|u^{\epsilon_m}|^2)|^2u^{\epsilon_m}- |g(|u^{\epsilon_n}|^2)|^2u^{\epsilon_n}\Big)dt\\
		&\quad-\bi \sum_{k\in \mathbb N^+} Im(Q^{\frac 12}e_k) Q^{\frac 12}e_k \Big(g'(|u^{\epsilon_m}|^2)g(|u^{\epsilon_m}|^2)|u^{\epsilon_m}|^2u^{\epsilon_m}- g'(|u^{\epsilon_n}|^2)g(|u^{\epsilon_n}|^2)|u^{\epsilon_n}|^2u^{\epsilon_n}\Big)dt\\
		&\quad+\mathbf i \Big(g(|u^{\epsilon_m}|^2)u^{\epsilon_m}-g(|u^{\epsilon_n}|^2)u^{\epsilon_n}\Big)dW(t)
		\end{align*}
		for multiplicative noise case. 
		Then using the It\^o formula to $\|u^{\epsilon_m}(t)-u^{\epsilon_n}(t)\|^2$, the properties of $f_{\epsilon}$ in Lemma \ref{sec-reg}, the mean value theorem and the Gagliardo--Nirenberg interpolation inequality, we obtain that for $\eta'(2-d)\le 2,$
		\begin{align}\label{main-add}
		&\quad\|u^{\epsilon_m}(t)-u^{\epsilon_n}(t)\|^2\\\nonumber
		&=\int_0^t 2\<u^{\epsilon_m}-u^{\epsilon_n},\bi \lambda f_{\epsilon_m}(|u^{\epsilon_m}|^2)u^{\epsilon_m}-f_{\epsilon_n}(|u^{\epsilon_n}|^2)u^{\epsilon_n}\>ds\\\nonumber
		&\le \int_0^t  4|\lambda| \|u^{\epsilon_m}(s)-u^{\epsilon_n}(s)\|^2 ds
		+4|\lambda| \int_0^t |Im \<u^{\epsilon_m}(s)-u^{\epsilon_n}(s),(f_{\epsilon_m}(|u^{\epsilon_n}|^2)-f_{\epsilon_n}(|u^{\epsilon_n}|^2)u^{\epsilon_n}\>| ds\\\nonumber
		&\le \int_0^t  6|\lambda| \|u^{\epsilon_m}(s)-u^{\epsilon_n}(s)\|^2 ds+4|\lambda| \int_0^t \|u^{\epsilon_m}(s)-u^{\epsilon_n}(s)\|_{L^1}
		\Big\|\frac {(\epsilon_m-\epsilon_n)|u^{\epsilon_n}|}{\epsilon_m+|u^{\epsilon_n}|^2}\Big\|_{L^{\infty}} ds\\\nonumber
		&\quad +2|\lambda|\int_0^t \|\log(1+\frac {(\epsilon_n-\epsilon_m)|u^{\epsilon_n}|^2}{1+\epsilon_m|u^{\epsilon_n}|^2}) |u^{\epsilon_n}|\|^2ds  \\\nonumber
		&\le \int_0^t  6|\lambda| \|u^{\epsilon_m}(s)-u^{\epsilon_n}(s)\|^2 ds
		+4|\lambda|\epsilon_n^{\frac 12} \int_0^t \|u^{\epsilon_m}(s)-u^{\epsilon_n}(s)\|_{L^1}d+2|\lambda| C\epsilon_n^{\eta'}\int_0^t\|u^{\epsilon_n}\|_{L^{2+2\eta'}}^{2+2\eta'}ds
		\end{align}
		for additive noise case,
		and 
		\begin{align*}
		&\quad\|u^{\epsilon_m}(t)-u^{\epsilon_n}(t)\|^2\\
		&=\int_0^t 2\<u^{\epsilon_m}-u^{\epsilon_n},\bi \lambda f_{\epsilon_m}(|u^{\epsilon_m}|^2)u^{\epsilon_m}-f_{\epsilon_n}(|u^{\epsilon_n}|^2)u^{\epsilon_n}\>ds\\
		&\quad-\int_0^t \<u^{\epsilon_m}-u^{\epsilon_n},\sum_{k\in\mathbb N^+}|Q^{\frac 12}e_k|^2\Big(|g(|u^{\epsilon_m}|^2)|^2u^{\epsilon_m}- |g(|u^{\epsilon_n}|^2)|^2u^{\epsilon_n}\Big)\>ds\\
		&\quad-2\int_0^t\<u^{\epsilon_m}-u^{\epsilon_n},\bi\sum_{k\in \mathbb N^+}Im(Q^{\frac 12}e_k) Q^{\frac 12}e_k \Big(g'(|u^{\epsilon_m}|^2)g(|u^{\epsilon_m}|^2)|u^{\epsilon_m}|^2u^{\epsilon_m}- g'(|u^{\epsilon_n}|^2)g(|u^{\epsilon_n}|^2)|u^{\epsilon_n}|^2u^{\epsilon_n}\Big)\>ds\\
		&\quad+2\int_0^t\<u^{\epsilon_m}-u^{\epsilon_n}, \mathbf i \Big(g(|u^{\epsilon_m}|^2)u^{\epsilon_m}-g(|u^{\epsilon_n}|^2)u^{\epsilon_n}\Big)dW(s)\>\\
		&\quad+\int_0^t \<g(|u^{\epsilon_m}|^2)u^{\epsilon_m}- g(|u^{\epsilon_n}|^2)u^{\epsilon_n},\sum_{k\in\mathbb N^+}|Q^{\frac 12}e_k|^2\Big(g(|u^{\epsilon_m}|^2)u^{\epsilon_m}- g(|u^{\epsilon_n}|^2)u^{\epsilon_n}\Big)\>ds\\
		&\le \int_0^t  4|\lambda| \|u^{\epsilon_m}(s)-u^{\epsilon_n}(s)\|^2 ds
		+4|\lambda|\epsilon_n^{\frac 12} \int_0^t \|u^{\epsilon_m}(s)-u^{\epsilon_n}(s)\|_{L^1}ds +2|\lambda| C\epsilon_n^{\eta'}\int_0^t \|u^{\epsilon_n}\|_{L^{2+2\eta'}}^{2+2\eta'}ds\\
		&\quad+2\int_0^t\<u^{\epsilon_m}-u^{\epsilon_n}, \mathbf i \Big(g(|u^{\epsilon_m}|^2)u^{\epsilon_m}-g(|u^{\epsilon_n}|^2)u^{\epsilon_n}\Big)dW(s)\>\\
		&\quad+\int_0^t \<(g(|u^{\epsilon_m}|^2)-g(|u^{\epsilon_n}|^2))u^{\epsilon_n},\sum_{k\in\mathbb N^+}|Q^{\frac 12}e_k|^2(g(|u^{\epsilon_m}|^2)- g(|u^{\epsilon_n}|^2))u^{\epsilon_m}\>ds\\
		&-2\int_0^t\<u^{\epsilon_m}-u^{\epsilon_n},\bi \sum_{k\in \mathbb N^+}Im(Q^{\frac 12}e_k) Q^{\frac 12}e_k\Big(g'(|u^{\epsilon_m}|^2)g(|u^{\epsilon_m}|^2)|u^{\epsilon_m}|^2u^{\epsilon_m}- g'(|u^{\epsilon_n}|^2)g(|u^{\epsilon_n}|^2)|u^{\epsilon_n}|^2u^{\epsilon_n}\Big)\>ds
		\end{align*}
		for multiplicative noise case. 	
		By  using \eqref{con-g} and \eqref{con-g1} in Assumption \ref{ass-strong} and the assumptions on $Q$, we have that 
		\begin{align}\label{main-mul}
		&\quad\|u^{\epsilon_m}(t)-u^{\epsilon_n}(t)\|^2\\\nonumber
		&\le \int_0^t  (4|\lambda|+C(g,Q)) \|u^{\epsilon_m}(s)-u^{\epsilon_n}(s)\|^2 ds
		+4|\lambda|\epsilon_n^{\frac 12} \int_0^t \|u^{\epsilon_m}(s)-u^{\epsilon_n}(s)\|_{L^1}ds\\\nonumber
		&\quad+2|\lambda| C\epsilon_n^{\eta'}\int_0^t \|u^{\epsilon_n}\|_{L^{2+2\eta'}}^{2+2\eta'}ds+\int_0^t\<u^{\epsilon_m}-u^{\epsilon_n}, \mathbf i \Big(g(|u^{\epsilon_m}|^2)u^{\epsilon_m}-g(|u^{\epsilon_n}|^2)u^{\epsilon_n}\Big)dW(s)\>.
		\end{align}
		
		Next we show the strong convergence of the sequence $\{u^{\epsilon_n}\}_{n\in\mathbb N^+}$ in the following different cases.
		
		Case 1: $\mathcal O$ is a bounded domain.
		By using the H\"older inequality $\|u^{\epsilon_m}(s)-u^{\epsilon_n}(s)\|_{L^1}\le |\mathcal O|^{\frac 12}\|u^{\epsilon_m}(s)-u^{\epsilon_n}(s)\|$ on \eqref{main-add} and \eqref{main-mul}, and using the Gronwall's inequality, we get  
		\begin{align*}
		\sup_{t\in[0,T]}\|u^{\epsilon_m}(t)-u^{\epsilon_n}(t)\|^2\le C(\lambda,T,|\mathcal O|) (\epsilon_n+\epsilon_n^{\eta'})(1+\sup_{t\in[0,T]}\|u^{\epsilon_n}\|_{L^{2+2\eta'}}^{2+2\eta'})
		\end{align*}
		for additive noise case. 
		In the multiplicative noise case, taking supreme over $t$ and then taking expectation on \eqref{main-mul}, together with the Burkholder and Young inequalities, we get that for a small $\kappa>0$,
		\begin{align*}
		&\quad\E\Big[\sup_{t\in [0,T]}\|u^{\epsilon_m}(t)-u^{\epsilon_n}(t)\|^2\Big]\\
		&\le C(\lambda,T,|\mathcal O|,Q)  (\epsilon_n+\epsilon_n^{\eta'})+C\E\Big[\sup_{t\in [0,T]}\Big|\int_0^t\<u^{\epsilon_m}-u^{\epsilon_n}, \mathbf i \Big(g(|u^{\epsilon_m}|^2)u^{\epsilon_m}-g(|u^{\epsilon_n}|^2)u^{\epsilon_n}\Big)dW(s)\>\Big|\Big]\\
		&\le C(\lambda,T,|\mathcal O|,Q)  (\epsilon_n+\epsilon_n^{\eta'})+C\E\Big[\Big(\int_0^T\sum_{i}\|Q^{\frac 12}e_i\|_{L^{\infty}}^2
		\|u^{\epsilon_m}-u^{\epsilon_n}\|^4ds\Big)^{\frac 12}\Big]\\
		&\le  C(\lambda,T,|\mathcal O|,Q)  (\epsilon_n+\epsilon_n^{\eta'})+C\E\Big[\sup_{s\in [0,T]}\|u^{\epsilon_m}-u^{\epsilon_n}\|\Big(\int_0^T\sum_{i}\|Q^{\frac 12}e_i\|_{L^{\infty}}^2
		\|u^{\epsilon_m}-u^{\epsilon_n}\|^2ds\Big)^{\frac 12}\Big]\\
		&\le C(\lambda,T,|\mathcal O|,Q)  (\epsilon_n+\epsilon_n^{\eta'})+\kappa\E\Big[\sup_{t\in [0,T]}\|u^{\epsilon_m}(t)-u^{\epsilon_n}(t)\|^2\Big]\\
		&\quad+C(\kappa)\E\Big[\int_0^T\sum_{i}\|Q^{\frac 12}e_i\|_{L^{\infty}}^2
		\|u^{\epsilon_m}-u^{\epsilon_n}\|^2ds\Big].
		\end{align*}
		Taking  $\kappa<\frac 12$, we have that 
		$$\E\Big[\sup_{t\in [0,T]}\|u^{\epsilon_m}(t)-u^{\epsilon_n}(t)\|^2\Big]\le C(Q,T,\lambda,p,u_0,|\mathcal O|)(\epsilon_n+\epsilon_n^{\eta'}).$$ 
		
		Case 2: $\mathcal O=\mathbb R^d$.
		Since $u_0\in L_2^{\alpha}, \alpha\in (0,1],$
		using the interpolation inequality in Lemma \ref{wei-hol}  implies that for any $\eta\in [0,1)$ and $\alpha>\frac {\eta d}{2(1-\eta)}$ (i.e., $\eta\in (0,\frac {2\alpha} {2\alpha+d})$),
		\begin{align*}
		&\|u^{\epsilon_m}(t)-u^{\epsilon_n}(t)\|^2\\\nonumber
		\le& \int_0^t  4|\lambda| \|u^{\epsilon_m}(s)-u^{\epsilon_n}(s)\|^2 ds
		+4|\lambda|\int_0^t \epsilon_{n}^{\frac \eta2}\|u^{\epsilon_m}(s)-u^{\epsilon_n}(s)\|\|u^{\epsilon_n}\|_{L^{2-2\eta}}^{1-\eta}ds\\\nonumber
		&+2|\lambda| C\epsilon_n^{\eta'}\int_0^t \|u^{\epsilon_n}\|_{L^{2+2\eta'}}^{2+2\eta'}ds\\\nonumber
		\le&  \int_0^t  4|\lambda| \|u^{\epsilon_m}(s)-u^{\epsilon_n}(s)\|^{2} ds
		+2|\lambda|\int_0^t\|u^{\epsilon_m}(s)-u^{\epsilon_n}(s)\|^2+ 
		2|\lambda| \epsilon_{m}^{\eta} \int_0^t\|u^{\epsilon_n}\|_{L^{2-2\eta}}^{2-2\eta}ds\\\nonumber
		&+2|\lambda| C\epsilon_n^{\eta'}\int_0^t \|u^{\epsilon_n}\|_{L^{2+2\eta'}}^{2+2\eta'}ds\\\nonumber
		\le&\int_0^t  6|\lambda| \|u^{\epsilon_m}(s)-u^{\epsilon_n}(s)\|^2 ds
		+2|\lambda|C \epsilon_{n}^{\eta} \int_0^t\|u^{\epsilon_n}\|_{L^{2}_\alpha}^{\frac {d\eta}{\alpha}}\|u^{\epsilon_n}\|^{2-2\eta-\frac{d\eta}{\alpha} }ds\\\nonumber
		&+2|\lambda| C\epsilon_n^{\eta'}\int_0^t \|u^{\epsilon_n}\|_{L^{2+2\eta'}}^{2+2\eta'}ds
		\end{align*}
		for additive noise case,
		and 
		\begin{align*}
		&\|u^{\epsilon_m}(t)-u^{\epsilon_n}(t)\|^2\\\nonumber
		\le&\int_0^t  C\|u^{\epsilon_m}(s)-u^{\epsilon_n}(s)\|^2 ds
		+C \epsilon_{m}^{\eta} \int_0^t\|u^{\epsilon_n}\|_{L^{2}_\alpha}^{\frac {d\eta}{\alpha}}\|u^{\epsilon_n}\|^{2-2\eta-\frac{d\eta}{\alpha} }ds\\\nonumber
		&+2|\lambda| C\epsilon_n^{\eta'}\int_0^t \|u^{\epsilon_n}\|_{L^{2+2\eta'}}^{2+2\eta'}ds+\Big|\int_0^t\<u^{\epsilon_m}-u^{\epsilon_n}, \mathbf i \Big(g(|u^{\epsilon_m}|^2)u^{\epsilon_m}-g(|u^{\epsilon_n}|^2)u^{\epsilon_n}\Big)dW(s)\>\Big|\nonumber
		\end{align*}
		for multiplicative noise case.
		Then taking supreme  over $t$, taking expectation, using \eqref{con-g}, Lemma \ref{mass-pri} and Proposition \ref{wei-hol-1}, and applying Gronwall's inequality, we have that for $\alpha\in (0,1],$ $\eta\in (0,\frac {2\alpha} {2\alpha+d})$ and $\eta'(d-2)\le 2$,
		\begin{align}\label{main-rd-sec-wei}
		&\E\Big[\sup_{t\in [0,T]}\|u^{\epsilon_m}(t)-u^{\epsilon_n}(t)\|^2\Big]\\\nonumber
		\le& C(T,Q,u_0,g)\E\Big[\sup_{[0,T]}\Big(\|u^{\epsilon_n}\|_{L^{2}_\alpha}^{\frac {d\eta}{\alpha}}\|u^{\epsilon_n}\|^{2-2\eta-\frac{d\eta}{\alpha} }+\|u^{\epsilon_n}\|_{L^{2+2\eta'}}^{2+2\eta'}\Big)\Big](\epsilon_n^{\eta}+\epsilon_n^{\eta'})\\\nonumber
		\le& C(T,Q,u_0,g,\alpha,\eta)\epsilon_n^{\min(\eta,\eta')}.
		\end{align}

		Step 2:  The limit process $u$ of $\{u^{\epsilon_n}\}_{n\in \mathbb N^+}$ in $\mathbb M_{\mathcal F}^2(C([0,T];\mathbb H))$ satisfies \eqref{SlogS} in mild form.
		We use the multiplicative noise case to present all the detailed procedures. 
		It suffices to prove that each term in the mild form of RSlogS equation \eqref{Reg-SlogS}
		convergenes to the corresponding part in 
		\begin{align*}
		&S(t)u_0+\bi\lambda\int_0^tS(t-s)\log(|u|^2)uds-\frac 12\int_0^tS(t-s)(g(|u|^2))^2u
		\sum_{k}|Q^{\frac 12}e_k|^2ds\\\nonumber
		&-\bi \int_0^tS(t-s)g'(|u|^2)g(|u|^2)|u|^2u
		\sum_{k}Im(Q^{\frac 12}e_k)Q^{\frac 12}e_kds+\bi\int_0^tS(t-s)g(|u|^2)udW(s)\\
		&:=S(t)u_0+V_1+V_2+V_3+V_4.
		\end{align*}
		We first claim that all the terms $V_1$-$V_4$ make sense. 
		By Lemma \ref{h1-pri} and Proposition \ref{wei-hol-1}, we have that for $p\ge2,$
		\begin{align*}
		\sup_{n}\E\left[\sup_{t\in[0,T]}\|u^{\epsilon_n}(t)\|_{\mathbb H^1}^p\right]+\sup_{n}\E\left[\sup_{t\in[0,T]}\|u^{\epsilon_n}(t)\|_{L^2_{\alpha}}^p\right]&\le C(u_0,T,Q).
		\end{align*}
		By applying the Fourier transform and 
		Parseval's theorem, using the Fatou theorem and strong convergence of $(u^{\epsilon_n})_{n\in \mathbb N^+}$ in $\mathbb M_{\mathcal F}^2(C([0,T];\mathbb H))$, we obtain   
		\begin{align*}	
		&\E\Big[\sup_{t\in[0,T]}\|u(t)\|_{\mathbb H^1}^2\Big]+\E\Big[\sup_{t\in[0,T]}\|u(t)\|_{L^2_{\alpha}}^2\Big]\\
		\le& \sup_{n}\E\Big[\sup_{t\in[0,T]}\|u^{\epsilon_n}(t)\|_{\mathbb H^1}^2\Big]
		+\sup_{n}\E[\sup_{t\in[0,T]}\|u^{\epsilon_n}(t)\|_{L^2_{\alpha}}^2]\le C(u_0,T,Q).
		\end{align*}
		Then 
		the Gagliardo--Nirenberg interpolation inequality yields that for  small $\eta',\eta>0$, 
		\begin{align*}
		\|\log(|u|^2)u\|^2&=
		\int_{|u|^2\ge 1} (\log(|u|^2))^2|u|^2dx+\int_{|u|^2\le 1}(\log(|u|^2))^2|u|^2dx\\
		&\le \int_{|u|^2\ge 1} |u|^{2+{2\eta'}}dx+\int_{|u|^2\le 1}|u|^{2-{2\eta}}dx\\
		&\le C(\|u\|_{L^{2+2\eta}}^{2+2\eta'}+\|u\|_{L^{2-2\eta}}^{2-2\eta})\\
		&\le C(\|u\|^{2-2\eta}_{L^{2-2\eta}}+\|\nabla u\|^{d\eta'}\|u\|^{2\eta'+2-d\eta'}).
		\end{align*}
		When $\mathcal O=\mathbb R^d$, we use the weighted version of the interpolation inequality in Lemma \ref{wei-hol} to deal with the term $\|u\|_{L^{2-2\eta}}^{2-2\eta}$, 
		and have  that for small $\eta<\frac {2\alpha}{2\alpha+d}$,
		$$\|\log(|u|^2)u\|^2 \le C(\|\nabla u\|^{d\eta'}\|u\|^{2\eta'+2-d\eta'}+\|u\|_{L^2_{\alpha}}^{\frac {2d\eta}{\alpha}}\|u\|^{2-2\eta-\frac {2d\eta}{\alpha}}).$$
		This implies that $V_1$ makes sense in $\mathbb M_{\mathcal F}^2(C([0,T];\mathbb H))$ by Proposition \ref{wei-hol-1}, Lemmas \ref{mass-pri} and \ref{h1-pri}.
		Meanwhile, we can show that  $V_2$-$V_4 \in \mathbb M_{\mathcal F}^2(C([0,T];\mathbb H))$ by using the Minkowski and Burkerholder inequalities due to our assumption on $g$ and $Q$. 
		
		Next, we show that the mild form of $u^{\epsilon_n}$ converges to 
		$S(t)u_0+V_1+V_2+V_3+V_4.$ To prove that 
		\begin{align*}
		&\lim_{n\to\infty}\E\Big[\sup_{t\in[0,T]}\Big\|\int_0^tS(t-s)f_{\epsilon_n}(|u^{\epsilon_n}|^2)u^{\epsilon_n}ds-V_1\Big\|^2\Big]=0,
		\end{align*} 
		we use the following decomposition of $f_{\epsilon_n}(|u^{\epsilon_n}|^2)u^{\epsilon_n}-\log(|u|^2)u.$
		When $|u|>|u^{\epsilon_n}|$,  
		\begin{align*}
		&f_{\epsilon_n}(|u^{\epsilon_n}|^2)u^{\epsilon_n}-\log(|u|^2)u\\
		=&(f_{\epsilon_n}(|u^{\epsilon_n}|^2)-f_{\epsilon_n}(|u|^2))u^{\epsilon_n}+f_{\epsilon_n}(|u|^2)(u^{\epsilon_n}-u)+(f_{\epsilon_n}(|u|^2)-\log(|u|^2))u,
		\end{align*}
		and when $|u|<|u^{\epsilon_n}|,$ 
		\begin{align*}
		&f_{\epsilon_n}(|u^{\epsilon_n}|^2)u^{\epsilon_n}-\log(|u|^2)u\\
		=&(\log(|u^{\epsilon_n}|^2)-\log(|u|^2))u+\log(|u^{\epsilon_n}|^2)(u^{\epsilon_n}-u)+(f_{\epsilon_n}(|u^{\epsilon_n}|^2)-\log(|u^{\epsilon_n}|^2))u^{\epsilon_n}.
		\end{align*}
		For convenience, let us show the estimate for $|u|>|u^{\epsilon_n}|$, the other case will be estimated in a similar way. 
		By using the H\"older inequality and the mean-valued theorem, we have that for small $\gamma>0,$
		\begin{align*}
		&\Big|(f_{\epsilon_n}(|u^{\epsilon_n}|^2)-f_{\epsilon_n}(|u|^2))u^{\epsilon_n}\Big|\\
		\le& \Big|(\frac {|u|^2-|u^{\epsilon_n}|^2}{\epsilon_n+|u^{\epsilon_n}|^2})^{\frac 12-\frac 12\gamma}(f_{\epsilon_n}(|u^{\epsilon_n}|^2)-f_{\epsilon_n}(|u|^2))^{\frac 12+\frac 12\gamma}u^{\epsilon_n}\Big|\\
		\le&  \Big|(|u|-|u^{\epsilon}_n|)^{\frac 12-\frac 12\gamma}(|u|+|u^{\epsilon}_n|)^{\frac 12-\frac 12\gamma} \frac {|u^{\epsilon_n}|}{(\epsilon_n+|u^{\epsilon_n}|^2)^{\frac 12-\frac 12\gamma}}
		(f_{\epsilon_n}(|u^{\epsilon_n}|^2)+f_{\epsilon_n}(|u|^2))^{\frac 12+\frac 12\gamma}\Big|.
		\end{align*}
		Then it implies that for  small enough $\eta>0,$
		\begin{align*}
		&\int_{|u|>|u^{\epsilon_n}|}\Big|(f_{\epsilon_n}(|u^{\epsilon_n}|^2)-f_{\epsilon_n}(|u|^2))u^{\epsilon_n}\Big|^2dx\\
		\le& \int_{|u|>|u^{\epsilon_n}|} \frac {|u^{\epsilon_n}|^2}{(\epsilon_n+|u^{\epsilon_n}|^2)^{1-\gamma}} |u-u^{\epsilon_n}|^{1-\gamma} 
		(|u|+|u^{\epsilon}|)^{1-\gamma}(f_{\epsilon_n}(|u^{\epsilon_n}|^2)+f_{\epsilon_n}(|u|^2))^{1+\gamma} dx\\
		\le& \int_{|u|>|u^{\epsilon_n}|, \epsilon_n+|u^{\epsilon_n}|\le 1} \frac {|u^{\epsilon_n}|^2}{(\epsilon_n+|u^{\epsilon_n}|^2)^{1-\gamma}} |u-u^{\epsilon_n}|^{1-\gamma} (|u|+|u^{\epsilon}|)^{1-\gamma}(f_{\epsilon_n}(|u^{\epsilon_n}|^2)+f_{\epsilon_n}(|u|^2))^{1+\gamma} dx\\
		& +\int_{|u|>|u|^{\epsilon_n}, \epsilon_n+|u^{\epsilon_n}|\ge 1} \frac {|u^{\epsilon_n}|^2}{(\epsilon_n+|u^{\epsilon_n}|^2)^{1-\gamma}} |u-u^{\epsilon_n}|^{1-\gamma} (|u|+|u^{\epsilon}|)^{1-\gamma}(f_{\epsilon_n}(|u^{\epsilon_n}|^2)+f_{\epsilon_n}(|u|^2))^{1+\gamma} dx\\
		\le&  C\int_{|u|>|u^{\epsilon_n}|, \epsilon_n+|u^{\epsilon_n}|\le 1} \frac {|u^{\epsilon_n}|^2}{(\epsilon_n+|u^{\epsilon_n}|^2)^{1-\gamma}} |u-u^{\epsilon_n}|^{1-\gamma} (|u|+|u^{\epsilon}|)^{1-\gamma}
		((\epsilon_n+|u^{\epsilon_n}|^2)^{-\eta}+(\epsilon_n+|u|^2)^{\eta})dx\\
		& +C\int_{|u|>|u|^{\epsilon_n}, \epsilon_n+|u^{\epsilon_n}|\ge 1} \frac {|u^{\epsilon_n}|^2}{(\epsilon_n+|u^{\epsilon_n}|^2)^{1-\gamma}} (|u|+|u^{\epsilon}|)^{1-\gamma} |u+u^{\epsilon}|^{1-\gamma}(\epsilon_n+|u|^2)^{\eta} dx.
		\end{align*}
		Now choosing $2-2\gamma+2\eta \le 2$, using the H\"older inequality and the weighted interpolation inequality in Lemma \ref{wei-hol}, we have that for $\alpha>\frac {\eta d}{1+\gamma-2\eta}$,
		\begin{align*}
		&\quad\int_{|u|>|u^{\epsilon_n}|}\Big|(f_{\epsilon_n}(|u^{\epsilon_n}|^2)-f_{\epsilon_n}(|u|^2))u^{\epsilon_n}\Big|^2dx\\
		&\le C\int_{\mathbb R^d}|u^{\epsilon_n}|^{2\gamma-2\eta} |u-u^{\epsilon_n}|^{1-\gamma} (|u|+|u^{\epsilon}|)^{1-\gamma} dx\\
		&\quad+C\int_{\mathbb R^d}|u^{\epsilon_n}|^{2\gamma-2\eta} |u-u^{\epsilon_n}|^{1-\gamma} (|u|+|u^{\epsilon}|)^{1-\gamma}(\epsilon_n+|u|^2)^{\eta} dx\\
		&\le C\|u-u^{\epsilon_n}\|^{1-\gamma}\Big\| (|u|^{2\gamma-2\eta}+|u|^{2\gamma})(|u|+|u^{\epsilon}|)^{1-\gamma}\Big\|_{L^{\frac 2{1+\gamma}}}\\
		&\le C\|u-u^{\epsilon_n}\|^{1-\gamma}\Big(\|u\|^{1+\gamma}+\|u\|_{L^{2-\frac {4\eta}{1+\gamma}}}^{1+\gamma-2\eta}\Big)\\
		&\le C\|u-u^{\epsilon_n}\|^{1-\gamma}\Big(\|u\|^{1+\gamma}+\|u\|_{L^{2}_{\alpha}}^{\frac {\eta d}{\alpha(1+r)}}\|u\|^{1+\gamma-2\eta-\frac {\eta d}{\alpha(1+\gamma)}}\Big).
		\end{align*}
		For the term $f_{\epsilon_n}(|u|^2)(u^{\epsilon_n}-u),$ we similarly have that for $\eta'<\frac {2}{\min(d-2,0)}$ and $\eta<\frac{2\alpha}{2\alpha+d},$ 
		\begin{align*}
		&\quad\int_{|u|>|u^{\epsilon}|}|f_{\epsilon_n}(|u|^2)(u^{\epsilon_n}-u)|^2dx\\
		&\le \int_{|u|>|u^{\epsilon}|,\epsilon_n+|u|^2\le 1} (\epsilon_n+|u|^2)^{-\frac \eta 2}|u^{\epsilon_n}-u|^2dx+\int_{|u|>|u^{\epsilon}|, \epsilon_n+|u|^2\ge 1} (\epsilon_n+|u|^2)^{\frac \eta 2}|u^{\epsilon_n}-u|^2dx\\
		&\le C\|u^{\epsilon_n}-u\| \Big(\|u^{\epsilon_n}-u\|+\|u\|_{L^{2+2\eta'}}^{1+\eta'} +\|u^{\epsilon_n}\|_{L^{2-2\eta}}^{1-\eta}+\|u\|_{L^{2-2\eta}}^{1-\eta}\Big)\\
		&\le C\|u^{\epsilon_n}-u\| \Big(\|u^{\epsilon_n}-u\|+\|\nabla u\|^{\frac {\eta' d}{2}}\|u\|^{1+\eta'-\frac {\eta' d}{2}}+\|u\|_{L^2_{\alpha}}^{\frac {d\eta}{2\alpha}}\|u\|^{1-\eta-\frac {d\eta}{2\alpha}}+\|u^{\epsilon_n}\|_{L^2_{\alpha}}^{\frac {d\eta}{2\alpha}}\|u^{\epsilon_n}\|^{1-\eta-\frac {d\eta}{2\alpha}}\Big).
		\end{align*}
		For the term $(f_{\epsilon_n}(|u|^2)-\log(|u|^2))u,$ the mean-valued theorem, the property that $\log(1+|x|)\le |x|$ and the Gagliardo--Nirenberg interpolation inequality
		\begin{align*}
		\|u\|_{L^q(\mathbb R^d)}\le C_1\|u\|^{1-\alpha}\|\nabla u\|^{\alpha}, \quad q=\frac {2d}{d-2\alpha},\; \text{for}  \;\alpha\in (0,1]
		\end{align*}
		with $q=2\eta'+2$, yield that for  $\eta'(d-2)\le 2$ and $\eta\le \frac{\alpha}{2\alpha-d},$
		\begin{align*}
		&\int_{\mathcal O}|(f_{\epsilon_n}(|u|^2)-\log(|u|^2))u|^2 dx\le C\epsilon_n^{\eta'}\|u\|_{L^{2\eta'+2}}^{2\eta'+2}+C\epsilon_n^{\eta}\|u\|_{L^{2-2\eta}}^{2-2\eta}\\
		\le& 
		C\epsilon_n^{\eta'}\|\nabla u\|^{\frac {\eta' d}{2}}\|u\|^{1+\eta'-\frac {\eta' d}{2}}+C\epsilon_n^{\eta}
		\|u^{\epsilon_n}\|_{L^2_{\alpha}}^{\frac {d\eta}{\alpha}}\|u^{\epsilon_n}\|^{1-\eta-\frac {d\eta}{\alpha}}.
		\end{align*} 
		Combining the above estimates, using the a priori estimate of $u^{\epsilon}$ and $u$ in Lemmas \ref{mass-pri} and \ref{h1-pri} and Proposition \ref{wei-hol-1}, and applying the strong convergence of $u^{\epsilon_n}$ \eqref{main-rd-sec-wei},we obtain that 
		\begin{align*}
		\lim_{n\to\infty}\E\Big[\sup_{t\in[0,T]}\Big\|\int_0^tS(t-s)f_{\epsilon_n}(|u^{\epsilon_n}|^2)u^{\epsilon_n}-\log(|u|^2)uds\Big\|^2\Big]=0.
		\end{align*} 
		
		The Minkowski inequality, \eqref{con-g} and \eqref{con-g1} yield that 
		\begin{align*}
		&\quad\Big\|-\int_0^tS(t-s)\frac 12(g(|u^{\epsilon_n}|^2))^2u^{\epsilon_n}
		\sum_{k}|Q^{\frac 12}e_k|^2ds-V_2\Big\|\\
		&\le \sum_{k}\|Q^{\frac 12}e_k\|_{L^{\infty}}^2\int_0^T\|g(|u^{\epsilon_n}|^2))^2u^{\epsilon_n}-g(|u|^2))^2u\|ds\\
		&\le C_gT\sum_{k}\|Q^{\frac 12}e_k\|_{L^{\infty}}^2\sup_{t\in[0,T]}\|u^{\epsilon_n}(t)-u(t)\|,
		\end{align*}
		and 
		\begin{align*}
		&\quad\Big\|-\bi \int_0^tS(t-s)(g'(|u^{\epsilon_n}|^2)g(|u^{\epsilon_n}|^2)|u^{\epsilon_n}|^2 u^{\epsilon_n}
		\sum_{k}Im (Q^{\frac 12}e_k)Q^{\frac 12}e_kds-V_3\Big\|\\
		&\le C_gT\sum_{k}\|Q^{\frac 12}e_k\|_{L^{\infty}}^2\sup_{t\in[0,T]}\|u^{\epsilon_n}(t)-u(t)\|.
		\end{align*}
		The Burkerholder inequality and the unitary property of $S(\cdot)$ yield that 
		\begin{align*}
		&\quad\E\Big[\sup_{t\in[0,T]}\Big\|\bi \int_0^tS(t-s)g(|u^{\epsilon_n}|^2)u^{\epsilon_n}dW(s)-V_4\Big\|^2\Big]\\
		&\le C\E\Big[ \int_0^T\sum_{k}\|Q^{\frac 12}e_k\|_{L^{\infty}}^2\|g(|u^{\epsilon_n}|^2)u^{\epsilon_n}-g(|u|^2)u\|^2ds\Big]
		\le  C\E\Big[\sup_{t\in[0,T]} \|u^{\epsilon_n}(t)-u(t)\|^2\Big].
		\end{align*}
		Combining the above estimates and the strong convergence of $u^{\epsilon_n}$, we complete the proof of step 2.
		
		Step 3: $u$ is independent of the choice of the sequence of $\{u^{\epsilon_n}\}_{n\in \mathbb N^+}.$ Assume that $\widetilde u$ and $u$ are two different limit processes of two different sequences of $\{u^{\epsilon_n}\}_{n\in\mathbb N^+}$ and $\{u^{\epsilon_m}\}_{m\in\mathbb N^+}$, respectively.
		Then by step 2, they both satisfies Eq. \eqref{SlogS}. By repeating the procedures in step 1, it is not hard to obtain that $\widetilde u=u.$
	\qed
	
	The procedures in the above proof immediately yield the following convergence rate result for $u^{\epsilon}$ in the regularized problem \eqref{Reg-SlogS} and the H\"older regularity estimate of $u^{\epsilon}$ and $u^0.$
	
	\begin{cor}\label{strong-con}
		Let the condition of Theorem \ref{mild-general} hold.  Assume that  $u^{\epsilon}$ is the mild solution in Proposition \ref{prop-well-reg}, $\epsilon\in (0,1).$ For $p\ge 2$, there exists $C(Q,T,\lambda,p,u_0)>0$ such that for any $\eta'(d-2)\le 2,$
		\begin{align*}
		\E\Big[\sup_{t\in [0,T]}\|u(t)-u^{\epsilon}(t)\|^p\Big] \le C(Q,T,\lambda, p,u_0)(\epsilon^{\frac p2}+\epsilon^{\frac {\eta'p} 2})
		\end{align*}
		when $\mathcal O$ is bounded domain, and
		\begin{align*}
		\E\Big[\sup_{t\in [0,T]}\|u(t)-u^{\epsilon}(t)\|^{p}\Big]
		&\le C(Q,T,\lambda, p,u_0,\alpha)(\epsilon^{\frac {\alpha p}{2\alpha+d}}+\epsilon^{\frac {\eta'p} 2})
		\end{align*}
		when $\mathcal O=\mathbb R^d$.
	\end{cor}

	\begin{cor}\label{stro-contin}
		Let the condition of Theorem \ref{mild-general} hold.  Assume that  $u^{\epsilon}$ is the mild solution in Proposition \ref{prop-well-reg}, $\epsilon\in (0,1)$ and $u^{0}$ is the mild solution of Eq. \eqref{SlogS}. For $p\ge 2$, there exists $C(Q,T,\lambda,p,u_0)>0$ such that for $\epsilon\in [0,1],$
		\begin{align*}
		\E\Big[\|u^{\epsilon}(t)-u^{\epsilon}(s)\|^p\Big] \le C(Q,T,\lambda, p,u_0)|t-s|^{\frac p2}.
		\end{align*}
	\end{cor}
	
	{\noindent\bf{Proof}} 
		By means of the mild form of $u^{\epsilon},$ $\epsilon\in [0,1)$, the priori estimates of $u^{\epsilon}$ in $\mathbb H^1\cap L_\alpha^2$ in Lemmas \ref{h1-pri} and \ref{wei-hol-1}, and in Step 2 of the proof of Theorem \ref{mild-general}, and the Burkholder inequalty, we obtain the desirable result.
	\qed
	
	\subsection{Well-posedness of SlogS equation with super-linearly
		growing diffusion coefficients}
	In this part, we extends the scope of $\widetilde g$, which allows the diffusion with super-linear growth, for the well-posedness of SlogS equation driven by conservative multiplicative noise. For instance, it includes the example $\widetilde g(x)=\bi g(|x|^2)x=\bi x\log(c+|x|^2),$ for $c>0$.
	
	\begin{tm}\label{mild-d1}
		Let  $W(t)$ be $L^2(\mathcal O;\mathbb R)$-valued and $g \in \mathcal C^1_b(\mathbb R)\cap \mathcal C(\mathbb R)$ satisfy the growth condition and  the embedding condition,
		\begin{align*}
		&\sup_{x\in [0,\infty)}|g'(x)x|\le C_g,\\
		&\|vg(|v|^2)\|\le C_d (1+\|v\|_{\mathbb H^1}+\|v\|_{L^2_{\alpha}})
		\end{align*}
		for some $q\ge 2$, $ \alpha\in [0,1]$, where $C_g>0$ depends on $g$, \;$C_d>0$ depends on $\mathcal O$, $d$, $\|v\|$ and $v\in \mathbb H^1\cap L^2_{\alpha}$. 
		Assume that $d=1$, $u_0\in \mathbb H^1\cap L^2_{\alpha}, $ $\alpha\in(0,1],$ and $\sum_{i\in \mathbb N^+}\|Q^{\frac 12}e_i\|_{\mathbb H^1}^2+\|Q^{\frac 12}e_i\|_{W^{1,\infty}}^2<\infty$. Then
		there exists a unique mild solution $u$ in $C([0,T];\mathbb H)$ for Eq. \eqref{SlogS} satisfying 
		\begin{align*}
		\E\Big[\sup_{t\in [0,T]}\|u(t)\|_{\mathbb H^1}^p\Big]+\E\Big[\sup_{t\in [0,T]}\|u(t)\|_{L^2_{\alpha}}^p\Big]\le C(Q,T,\lambda, p,u_0).
		\end{align*}
	\end{tm}
	
	{\noindent\bf{Proof}}
		By Proposition \ref{prop-well-reg} and Lemma \ref{h1-pri}, we can introduce the truncated sample space $$\Omega_R(t):=\left\{\omega: \sup_{s\in[0,t]}\|u^{\epsilon_m}(s)\|_{L^{\infty}}\le R, \sup_{s\in [0,t]}\|u^{\epsilon_n}(s)\|_{L^{\infty}}\le R\right\},$$
		where $n\le m.$
		The Gagliardo--Nirenberg interpolation inequality in $d=1$, the priori estimate in $\mathbb H^1$ and the continuity in $L^2$ of $u^{\epsilon_n}$ 
		imply that $u^{\epsilon_n}$ are continuous in $L^{\infty}$ a.s.
		Define a stopping time $$\tau_R:=\inf\{t\ge 0: \min(\sup_{s\in[0,t]}\|u^{\epsilon_m}(s)\|_{L^{\infty}},\sup_{s\in[0,t]}\|u^{\epsilon_n}(s)\|_{L^{\infty}})\ge R\}\wedge T.$$ Then on $\Omega_R(T),$ we have $\tau_R=T.$
		Let us take $f_{\epsilon}(x)=\log(x+\epsilon), x>0$ for convenience.
		It is obvious that $\Omega_R(t)\to \Omega$ as $R\to \infty$ and that for any $p\ge 1$,
		\begin{align*}
		&\quad\mathbb P\Big(\sup_{t\in[0,T]}\min(\|u^{\epsilon_m}(t)\|_{L^{\infty}},\|u^{\epsilon_n}(t)\|_{L^{\infty}})\ge R\Big)\\
		\le& C\frac 1{R^p}\Big(\E\Big[\sup_{t\in[0,T]}\|u^{\epsilon_m}(t)\|_{L^{\infty}}^p\Big]+\sup_{t\in[0,T]} \E\Big[\|u^{\epsilon_m}(t)\|_{L^{\infty}}^p\Big]\Big).
		\end{align*} 
		
		Step 1: $\{u^{\epsilon_n}\}_{n\in \mathbb N^+}$ forms a Cauchy sequence in  $ \mathbb M_{\mathcal F}^2(\Omega;C([0,T];\mathbb H^2))$. 
		Following the same steps like the proof of Theorem \ref{mild-general}, applying the It\^o formula on $\Omega_R(t)$ for $t\in (0,\tau_R)$ yields that 
		\begin{align*}
		&\quad\|u^{\epsilon_m}(t)-u^{\epsilon_n}(t)\|^2\\
		\le& \int_0^t  4|\lambda| \|u^{\epsilon_m}(s)-u^{\epsilon_n}(s)\|^2 ds
		+4|\lambda|\epsilon_n^{\frac 12} \int_0^t \|u^{\epsilon_m}(s)-u^{\epsilon_n}(s)\|_{L^1}ds\\
		&+\int_0^t \sum_{k\in\mathbb N^+}\<|Q^{\frac 12}e_k|^2 
		(g(|u^{\epsilon_n}|^2)-g(|u^{\epsilon_m}|^2))u^{\epsilon_m},(g(|u^{\epsilon_n}|^2)-g(|u^{\epsilon_m}|^2))u^{\epsilon_n}\>ds\\
		&+\int_0^t\<u^{\epsilon_m}-u^{\epsilon_n}, \mathbf i \Big(g(|u^{\epsilon_m}|^2)u^{\epsilon_m}-g(|u^{\epsilon_n}|^2)u^{\epsilon_n}\Big)dW(s)\>.
		\end{align*}
		Taking expectation on $\Omega_R(t)$ yields that 
		\begin{align*}
		&\quad\E\Big[\mathbb I_{\Omega_R(t)}\|u^{\epsilon_m}(t)-u^{\epsilon_n}(t)\|^2\Big]\\
		\le& \int_0^t  4|\lambda| \E\Big[\mathbb I_{\Omega_R(t)}\|u^{\epsilon_m}(s)-u^{\epsilon_n}(s)\|^2 \Big]ds
		+4|\lambda|\epsilon_n^{\frac 12} \int_0^t \E\Big[\mathbb I_{\Omega_R(t)} \|u^{\epsilon_m}(s)-u^{\epsilon_n}(s)\|_{L^1}\Big]ds\\
		&\quad+\int_0^t \sum_{k\in\mathbb N^+}|Q^{\frac 12}e_k|_{L^{\infty}}^2 \E\Big[\mathbb I_{\Omega_R(t)}\<
		(g(|u^{\epsilon_n}|^2)-g(|u^{\epsilon_m}|^2))u^{\epsilon_m},(g(|u^{\epsilon_n}|^2)-g(|u^{\epsilon_m}|^2))u^{\epsilon_n}\>\Big]ds.
		\end{align*}
		Making use of the assumptions on $g$, we get 
		\begin{align*}
		&\quad\E\Big[\mathbb I_{\Omega_R(t)}\|u^{\epsilon_m}(t)-u^{\epsilon_n}(t)\|^2\Big]\\
		&\le \int_0^t  4|\lambda| \E\Big[\mathbb I_{\Omega_R(t)}\|u^{\epsilon_m}(s)-u^{\epsilon_n}(s)\|^2 \Big]ds
		+4|\lambda|\epsilon_n^{\frac 12} \int_0^t \E\Big[\mathbb I_{\Omega_R(t)}\|u^{\epsilon_m}(s)-u^{\epsilon_n}(s)\|_{L^1}\Big]ds\\
		&\quad+C(Q)\int_0^t  \E\Big[\mathbb I_{\Omega_R(t)}\int_{\mathcal O} {(|u^{\epsilon_m}|^2+|u^{\epsilon_n}|^2)|u^{\epsilon_n}||u^{\epsilon_m}|} |u^{\epsilon_m}(t)-u^{\epsilon_n}(t)|^2dx\Big]ds\\
		&\le\int_0^t  4|\lambda| \E\Big[\mathbb I_{\Omega_R(t)} \|u^{\epsilon_m}(s)-u^{\epsilon_n}(s)\|^2 \Big]ds
		+4|\lambda|\epsilon_n^{\frac 12} \int_0^t \E\Big[\mathbb I_{\Omega_R(t)}\|u^{\epsilon_m}(s)-u^{\epsilon_n}(s)\|_{L^1}\Big]ds\\
		&\quad+C(Q)\int_0^t  \E\Big[(1+R^4)\mathbb I_{\Omega_R(t)}\|u^{\epsilon_m}(t)-u^{\epsilon_n}(t)\|^2\Big]ds.
		\end{align*}
		If $\mathcal O$ is bounded, then H\"older inequality and Gronwall's inequality yield that 
		\begin{align*}
		\E\Big[\mathbb I_{\Omega_R(t)} \|u^{\epsilon_m}(t)-u^{\epsilon_n}(t)\|^2\Big]
		&\le C(u_0,Q,T,\lambda)e^{(1+R^4 )T}\epsilon_n.
		\end{align*}
		On the other hand, the Chebyshev inequality and the a priori estimate lead to
		\begin{align*}
		&\E\Big[\mathbb I_{\Omega^c_R(t)} \|u^{\epsilon_m}(t)-u^{\epsilon_n}(t)\|^2\Big]
		\le (\mathbb P(\Omega^c_R(t)))^{\frac 1{p_1}}\Big(\E\Big[\|u^{\epsilon_m}(t)-u^{\epsilon_n}(t)\|^{2q_1}\Big]\Big)^{\frac 1q_1},
		\end{align*}
		where $\frac 1{p_1}+\frac 1{q_1}=1.$
		From the above estimate, choosing $p\gg p_1$ and denote by $\kappa=\frac {p}{p_1}$, we conclude that 
		\begin{align*}
		&\E\Big[ \|u^{\epsilon_m}(t)-u^{\epsilon_n}(t)\|^2\Big]\le C(u_0,Q,T,\lambda,p,p_1)\Big(e^{(1+R^4)T}\epsilon_n+R^{-\kappa}\Big)
		\end{align*}
		Then one may take $R=(\frac {c}{T} |\log(\epsilon)|)^{\frac 14}$ for $c\in (0,1)$ and get 
		\begin{align*}
		&\E\Big[ \|u^{\epsilon_m}(t)-u^{\epsilon_n}(t)\|^2\Big]\le  C(u_0,Q,T,\lambda,p,p_1)(\epsilon_n^{1-c}+(\frac {c}{T} |\log(\epsilon_n)|)^{-\frac \kappa 4}).
		\end{align*}
		By further applying the Burkerholder inequality to the stochastic integral, we achieve that  for any $\kappa>0,$
		\begin{align*}
		\E\Big[ \sup_{t\in [0,T]}\|u^{\epsilon_m}(t)-u^{\epsilon_n}(t)\|^2\Big]
		&\le  C(u_0,Q,T,\lambda,p,p_1)|\log(\epsilon_n)|^{-\frac \kappa 4}.
		\end{align*}
		When $\mathcal O=\mathbb R^d$, we just repeat the procedures in the proof of the case that $g$ is bounded and obtain that for $\eta\in (0,\frac{2\alpha}{2\alpha+d}]$ and $\alpha\in (0,1],$
		
		\begin{align*}
		&\quad\E\Big[\mathbb I_{\Omega_R(t)}\|u^{\epsilon_m}(t)-u^{\epsilon_n}(t)\|^2\Big]\\
		&\le\int_0^t  4|\lambda| \E\Big[\mathbb I_{\Omega_R(t)} \|u^{\epsilon_m}(s)-u^{\epsilon_n}(s)\|^2 \Big]ds
		+C \epsilon_{m}^{\eta} \int_0^t\E\Big[\mathbb I_{\Omega_R(t)} \|u^{\epsilon_n}\|_{L^{2}_\alpha}^{\frac {d\eta}{\alpha}}\|u^{\epsilon_n}\|^{2-2\eta-\frac{d\eta}{\alpha} }\Big]ds\\
		&\quad+C(Q)\int_0^t  \E\Big[(1+R^4)\mathbb I_{\Omega_R(t)}\|u^{\epsilon_m}(t)-u^{\epsilon_n}(t)\|^2\Big]ds.
		\end{align*}
		By using Gronwall's inequality and the estimate of $\mathbb P(\mathbb I_{\Omega_R(t)}^c)$, we immediately have that for $\eta\in (0,\frac {2\alpha}{2\alpha+d}), \alpha\in (0,1]$ and any $\kappa>0,$
		\begin{align*}
		&\E\Big[ \|u^{\epsilon_m}(t)-u^{\epsilon_n}(t)\|^2\Big]\le C(u_0,Q,T,\lambda,p,p_1)\Big(e^{(1+R^4)T}\epsilon_n^{\eta}+R^{-\kappa}\Big)
		\end{align*}
		Taking $R=(\frac {\eta c}{T} |\log(\epsilon_n)|)^{\frac 14}$ for $c\eta \in (0,1)$ and using the Burkerholder inequality, we have for $\eta\in (0,\frac {2\alpha}{2\alpha+d}), \alpha\in (0,1]$ and any $\kappa>0,$
		\begin{align*}
		&\E\Big[ \sup_{t\in[0,T]}\|u^{\epsilon_m}(t)-u^{\epsilon_n}(t)\|^2\Big]\le C(u_0,Q,T,\lambda,p,p_1,\eta)|\log(\epsilon_n)|^{-\frac {\kappa}4}.
		\end{align*}
		
		Step 2. $u$ is the mild solution.
		
		Let us use the same notations and procedures as in step 2 of the proof in Theorem \ref{mild-general}. 
		To show that the mild form $u^{\epsilon_n}$ converges to 
		$S(t)u_0+V_1+V_2+V_3+V_4.$ We only need to estimate $V_2$ and $V_4$ since $V_3=0$.
		Define 
		$$\Omega_{R_1}(t):=\{\omega: \sup_{s\in[0,t]}\|u^{\epsilon_n}(s)\|_{L^{\infty}}\le R_1, \sup_{s\in [0,t]}\|u(s)\|_{L^{\infty}}\le R_1\},$$
		and a stopping time $$\tau_{R_1}:=\inf\{t\ge 0: \min(\sup_{s\in[0,t]}\|u(s)\|_{L^{\infty}},\sup_{s\in[0,t]}\|u^{\epsilon_n}(s)\|_{L^{\infty}})\ge R_1\}\wedge T.$$ Then on $\Omega_{R_1}(T),$ we have $\tau_{R_1}=T$. The Minkowski inequality and the properties of $g$ yield that on $\Omega_{R_1}(t),$ for a small enough $\eta_1>0$,
		\begin{align*}
		&\Big\|-\frac 12\int_0^tS(t-s)(g(|u^{\epsilon_n}|^2))^2u^{\epsilon_n}
		\sum_{k}|Q^{\frac 12}e_k|^2ds-V_2\Big\|\\
		\le& \sum_{k}|Q^{\frac 12}e_k|_{L^{\infty}}^2\int_0^T\|g(|u^{\epsilon_n}|^2))^2u^{\epsilon_n}-g(|u|^2))^2u\|ds\\
		\le& C_gT\sum_{k}|Q^{\frac 12}e_k|_{L^{\infty}}^2\Big(1+R_1^2\Big)\sup_{t\in[0,T]}\|u^{\epsilon_n}-u\|.
		\end{align*}
		On the other hand, for any $p_2>0$,
		\begin{align*}
		\E\Big[\mathbb I_{\Omega_{R_1}^c(t)}\Big\|-\frac 12\int_0^tS(t-s)(g(|u^{\epsilon_n}|^2))^2u^{\epsilon_n}
		\sum_{k}|Q^{\frac 12}e_k|^2ds-V_2\Big\|^2\Big]\le C(u_0,Q,T,p_2)R_1^{-p_2}.
		\end{align*}
		Taking $R_1=\mathscr O(|\log(\epsilon_n)|^{\frac{ \kappa_1} {4(4+p_2)}}),\kappa_1<\kappa,$ we have that 
		\begin{align*}
		\lim_{n\to \infty}\sup_{t\in[0,T]}\E\Big[\Big\|\int_0^t-\frac 12S(t-s)(g(|u^{\epsilon_n}|^2))^2u^{\epsilon_n}
		\sum_{k}|Q^{\frac 12}e_k|^2ds-V_2\Big\|^2\Big]=0,
		\end{align*}
		which immediately implies that 
		\begin{align*}
		\lim_{n\to \infty}\E\Big[\sup_{t\in[0,T]}\Big\|-\frac 12\int_0^tS(t-s)(g(|u^{\epsilon_n}|^2))^2u^{\epsilon_n}
		\sum_{k}|Q^{\frac 12}e_k|^2ds-V_2\Big\|^2\Big]=0.
		\end{align*}
		The Burkerholder inequality and the unitary property of $S(\cdot)$ yield that 
		\begin{align*}
		&\E\Big[\sup_{t\in[0,\tau_R]}\Big\|\int_0^t\bi S(t-s)g(|u^{\epsilon_n}|^2)u^{\epsilon_n}dW(s)-V_4\Big\|^2\Big]\\
		\le& C\E\Big[ \int_0^T\sum_{k}\|Q^{\frac 12}e_k\|_{L^{\infty}}^2\|g(|u^{\epsilon_n}|^2)u^{\epsilon_n}-g(|u|^2)u\|^2ds\Big]\\
		\le&  C(1+R_1^{2})\E\Big[\sup_{t\in[0,T]} \|u^{\epsilon_n}(t)-u(t)\|^2\Big].
		\end{align*}
		On the other hand, the Chebyshev inequality, together with the a priori estimate of $u^{\epsilon_n}$, implies that 
		\begin{align*}
		&\E\Big[\sup_{t\ge \tau_{R_1}}\Big\|\int_0^t\bi S(t-s)g(|u^{\epsilon_n}|^2)u^{\epsilon_n}dW(s)-V_4\Big\|^2\Big]
		\\
		=&\E\Big[ \sup_{t\in [0,T]}\mathbb I_{\Omega_{R_1}^c} \Big\|\int_0^tS(t-s)g(|u^{\epsilon_n}|^2)u^{\epsilon_n}dW(s)-V_4\Big\|^2\Big]\le C(u_0,Q,T,p)R_1^{-p_2}.
		\end{align*}
		Taking $R_1=\mathscr O(|\log(\epsilon_n)|^{\frac{ \kappa_1} {4(2+p_2)}}),\kappa_1<\kappa,$ we have that 
		\begin{align*}
		\lim_{n\to \infty}\E\Big[\sup_{t\in[0,T]}\Big\|\int_0^t\bi S(t-s)g(|u^{\epsilon_n}|^2)u^{\epsilon_n}dW(s)-V_4\Big\|^2\Big]
		=0.
		\end{align*}
		Combining the above estimates and the strong convergence of $u^{\epsilon_n}$, we complete the proof.
	\qed

	\begin{rk}
		One may extend the scope of $\widetilde g$ to an abstract framework by similar arguments. Here the assumption $d=1$ lies on the fact that in $\mathbb H^1$ is an algebra by Sololev embedding theorem. When considering the case $d\ge2$, one may use  $\mathbb H^{\bs},{\bs>\frac d2}$ as the underlying space for the local well-posedness. However, as stated in Lemma \ref{h2-pri}, it seems impossible to get the uniform bound of $u^{\epsilon}$ in $\mathbb H^{\bs}$ for $\bs\ge 2.$
	\end{rk}

	\section{Appendix}
	
	The original problem and regularized problem can be rewritten into the equivalent evolution forms
	\begin{align}\label{evo}
	du&=A udt+F(u)dt+G(u)dW(t),\\\nonumber
	u(0)&=u_0,
	\end{align}
	where $A=\mathbf i \Delta,$ $F$ is the Nemystkii operator of drift coefficient function, and $G$ are the Nemystkii operator of diffusion coefficient function. Then the mild solution of the above evolution is defined as follows.
	
	\begin{df}\label{global}
		A continuous $\mathbb H$-valued $\mathcal F_t$ adapted process $u$ is a solution to \eqref{evo} if it satisfies $\mathbb P$-a.s for all $t\in [0,T],$
		\begin{align*}
		u(t)=S(t)u_0+\int_0^t S(t-s) F(u(s))ds+\int_0^tG(u(s))dW(s),
		\end{align*}
		where $S(t)$ is the $C_0$-group generated by $A.$
	\end{df}
	
	\begin{df}\label{local}
		A local mild solution of \eqref{evo} is $(u,\tau):=(u,\tau_n,\tau)$ satisfying $\tau_n\nearrow \tau, a.s.,$ as $n\to \infty$, $u \in \mathbb M_{\mathcal F}^p(\Omega;C([0,\tau);\mathbb H^{\bs}), \bs>0, p\ge 1$ and that 
		\begin{align*}
		u(t)&=S(t)u_0+\int_0^t S(t-s)F(u(s))ds\\
		&+\int_0^t S(t-s) G(u(s))dW(s), a.s.,
		\end{align*}
		for $t\le \tau_n$ in $\mathbb H^2$ for $n\in \mathbb N^+$.
		Solutions of \eqref{evo} are called unique, if $$\mathbb P\Big(u_1(t)=u_2(t), \forall t\in [0,\sigma_1\wedge \sigma_2)\Big)=0.$$
		for all local mild solution $(u_1,\sigma_1)$ and $(u_2,\sigma_2).$ The local solution $(u,\tau)$ is called a global mild solution if $\tau=T, a.s.$ and $u\in \mathbb M_{\mathcal F}^p(\Omega;C([0,T];\mathbb H^{\bs}).$
	\end{df}
	
	\begin{lm}\label{frist-loc}
		Let $\epsilon\in (0,1)$. Then $f_{\epsilon}(x)=\log({|x|^2+\epsilon}), x\in\mathbb C,$ satisfies 
		\begin{align*}
		|Im[(f_{\epsilon}(x_1)-f_{\epsilon}(x_2))(\bar x_1-\bar x_2)]|\le 4|x_1-x_2|^2,
		\end{align*}
	\end{lm}
	
	\textbf{Proof}
		Without loss of generality, we assume that $0<|x_2|\le |x_1|.$ 
		Notice that
		\begin{align*}
		Im[(f_{\epsilon}(x_1)-f_{\epsilon}(x_2))(\bar x_1-\bar x_2)]
		=\frac 12(\log(\epsilon+|x_1|^2)-\log(\epsilon+|x_2|^2))Im(\bar x_1x_2-\bar x_2x_1).
		\end{align*}
		Direct calculation yields that 
		\begin{align*}
		|Im(\bar x_1x_2-\bar x_2x_1)|\le 2|x_2||x_1-x_2|.
		\end{align*}
		Using the fact that 
		\begin{align*}
		|\log(\epsilon+|x_1|^2)-\log(\epsilon+|x_2|^2)|&=2|\log((\epsilon+|x_1|^2)^{\frac 12})-\log((\epsilon+|x_2|^2)^{\frac 12})|,
		\end{align*}
		we obtain 
		\begin{align*}
		&|Im[(f_{\epsilon}(x_1)-f_{\epsilon}(x_2))(\bar x_1-\bar x_2)]|\\
		&\le 2|\log((\epsilon+|x_1|^2)^{\frac 12})-\log((\epsilon+|x_2|^2)^{\frac 12})||x_2||x_1-x_2|.
		\end{align*}
		The mean value theorem leads to the desired result.
	\qed

	\begin{lm}\label{sec-reg}
		Let $\epsilon\in (0,1)$. Then $f_{\epsilon}(|x|^2)=\log(\frac{|x|^2+\epsilon}{1+\epsilon |x|^2})$ satisfies the following properties,
		\begin{align*}
		|f_{\epsilon}(|x|^2)|&\le  |\log(\epsilon)|,\\
		|d_{|x|}f_{\epsilon}(|x|^2)|&\le \frac {2(1-\epsilon^2)|x|}{(\epsilon+|x|^2)(1+\epsilon|x|^2)},\\
		|Im[(f_{\epsilon}(|x_1|^2)x_1-f_{\epsilon}(|x_2|^2)x_2)(\bar x_1-\bar x_2)]|&\le 4(1-\epsilon^2)|x_1-x_2|^2.
		\end{align*}
	\end{lm}
	\textbf{Proof}
		The proof of first and second estimates are derived by the property of $\log(\cdot).$ The last estimate is proven by similar arguments in the proof of  Lemma \ref{frist-loc}. 
	\qed

	\textbf{[Proof of Proposition \ref{eng-pri}]} Due to Lemma \ref{h1-pri}, it suffices to prove 
		\begin{align*}
		\E\Big[\sup\limits_{t\in[0,T]} (F_{\epsilon}(|u^{\epsilon}(t)|^2))^{p}\Big]\le C(u_0,T,Q,p).
		\end{align*}
		Let us take $f_{\epsilon}(|x|^2)=\log(|x|^2+\epsilon)$ as an example to illustrate the procedures. The desirable estimate in case that 
		$f_{\epsilon}(|x|^2)=\log(\frac{|x|^2+\epsilon}{1+\epsilon |x|^2})$ can be obtained similarly.
		Using the property of logarithmic function, we have that for small $\eta>0,$
		\begin{align*}
		|F_{\epsilon}(|u^{\epsilon}(t)|^2)|
		&=\Big|\int_{\mathcal O}\Big((\epsilon+|u^{\epsilon}(t)|^2)\log(\epsilon+|u^{\epsilon}(t)|^2)-|u^{\epsilon}(t)|^2-\epsilon\log(\epsilon)\Big)dx\Big|\\
		&\le \|u^{\epsilon}(t)\|^2+\int_{\mathcal O}|u^{\epsilon}(t)|^2\log(\epsilon+|u^{\epsilon}(t)|^2)dx+\Big|\int_{\mathcal O}\epsilon(\log(\epsilon+|u^{\epsilon}(t)|^2)-\log(\epsilon))dx\Big|\\
		&\le 2\|u^{\epsilon}(t)\|^2+\| u^{\epsilon}(t)\|_{L^{2-2\eta}}^{2-2\eta}+\| u^{\epsilon}(t)\|_{L^{2+2\eta}}^{2+2\eta},
		\end{align*}
		where we use the following estimation,
		for any small enough $\eta>0$,
		\begin{align*}
		&\int_{\mathcal O}|u^{\epsilon}(t)|^2\log(\epsilon+|u^{\epsilon}(t)|^2)dx\\
		&=
		\int_{\epsilon +|u^{\epsilon}|^2\ge 1} f_{\epsilon}(|u^{\epsilon}|^2)|u^{\epsilon}|^2dx+\int_{\epsilon +|u^{\epsilon}|^2\le 1} f_{\epsilon}(|u^{\epsilon}|^2)|u^{\epsilon}|^2dx\\
		&\le \int_{\epsilon +|u^{\epsilon}|^2\ge 1} (\epsilon +|u^{\epsilon}|^2)^{2\eta}|u^{\epsilon}|^2dx+\int_{\epsilon +|u^{\epsilon}|^2\le 1}(\epsilon+|u^{\epsilon}|^2)^{-2\eta}|u^{\epsilon}|^2dx\\
		&\le \|u\|_{L^{2-2\eta}}^{2-2\eta}+\|u\|_{L^{2+2\eta}}^{2+2\eta}.
		\end{align*}
		Then by 
		the Gagliardo--Nirenberg interpolation inequality in  a bounded domain $\mathcal O$, i.e.,  
		\begin{align*}
		\|u\|_{L^q(\mathcal O)}\le C_1\|u\|^{1-\alpha}\|\nabla u\|^{\alpha}+C_2\|u\|, \quad q=\frac {2d}{d-2\alpha}, \; \text{for} \; \alpha\in (0,1],
		\end{align*}
		we have that
		\begin{align*}
		\int_{\mathcal O}F_{\epsilon}(|u^{\epsilon}(s)|^2)dx&\le |\int_{\mathcal O}\Big((\epsilon+|u^{\epsilon}|^2)\log(\epsilon+|u^{\epsilon}|^2)-|u^{\epsilon}|^2-\epsilon\log\epsilon\Big) dx|\\
		&\le C(1+\|u^{\epsilon}\|_{L^{2+2\eta}}^{2+2\eta})\le C(\|u^{\epsilon}\|+\|\nabla u^{\epsilon}\|^{\frac {d\eta}{2\eta+2}}\|u^{\epsilon}\|^{1-\frac {d\eta}{2\eta+2}}).
		\end{align*}
		Taking $p$th moment and applying Lemma \ref{h1-pri}, we complete the proof for the case that $\mathcal O$ is a bounded domain.
		
		When $\mathcal O=\mathbb R^d,$ we need to control $\|u^{\epsilon}\|_{L^{2-2\eta}}.$ By using the weighted interpolation inequality in Lemma \ref{wei-hol}
		with $\alpha>\frac {d\eta}{2-2\eta}, \alpha\in (0,1]$,
		and applying the Gagliardo--Nirenberg interpolation inequality,
		\begin{align*}
		\|u\|_{L^q(\mathbb R^d)}\le C_1\|u\|^{1-\alpha}\|\nabla u\|^{\alpha}, \quad q=\frac {2d}{d-2\alpha},\; \text{for}  \;\alpha\in (0,1],
		\end{align*}
		where $q=2\eta+2.$
		Based on Lemma \ref{h1-pri} and Lemma \ref{wei-hol-1}, we complete the proof by using the Young inequality and  taking $p$th moment. 
	\qed

	\textbf{[Sketch Proof of Lemma \ref{h2-pri}]}
		Due to the loss of the regularity of the solution in time, we can not establish the bound in $\mathbb H^2$ through $\frac {\partial u^{\epsilon}} {\partial t}$ like in the deterministic case. 
		According to Lemma \ref{h1-pri}, it suffices to bound $\|\Delta u^{\epsilon}\|^2.$ We present the procedures of the estimation of $\E[\|\Delta u\|^2]$ for the conservative multiplicative noise case. One can easily follow the procedures to obtain the estimate of $\E[\sup_{t\in[0,\tau]}\|\Delta u(t)\|^2]$ for both additive and multiplicative noises.
		
		By using the It\^o formula to $\|\Delta u^{\epsilon}\|^2$ we obtain that
		\begin{align*}
		\|\Delta u^{\epsilon}(t)\|^2&= \|\Delta u^{\epsilon}_0\|^2
		+2\int_0^t \<\Delta u^{\epsilon}(s),II_{det} \>ds\\
		&\quad+ 2\int_0^t \<\Delta u^{\epsilon}(s),II_{mod} \>ds+2II_{Sto},
		\end{align*}
		where 
		\begin{align*}
		II_{Sto}&:=\int_0^t\<\Delta u^{\epsilon},\mathbf i g(|u^{\epsilon}|^2)\nabla u^{\epsilon}\nabla dW(s)\>+\int_0^t\<\Delta u^{\epsilon},\mathbf i g(|u^{\epsilon}|^2)u^{\epsilon}\Delta dW(s)\>\\
		&+\int_0^t\<\Delta u^{\epsilon},\mathbf i2g'(|u^{\epsilon}|^2)Re(\bar u^{\epsilon} \nabla u^{\epsilon})\nabla u^{\epsilon} dW(s)+\int_0^t\<\Delta u^{\epsilon}, \mathbf i2g'(|u^{\epsilon}|^2)Re(\bar u^{\epsilon} \nabla u^{\epsilon})u^{\epsilon}\nabla dW(s)\>\\
		&+\int_0^t\<\Delta u^{\epsilon},\mathbf i4g''(|u^{\epsilon}|^2)(Re(\bar u^{\epsilon} \nabla u^{\epsilon}))^2u^{\epsilon} dW(s)+\int_0^t\<\Delta u^{\epsilon},\mathbf i 2g'(|u^{\epsilon}|^2)(Re(\bar u^{\epsilon} \Delta u^{\epsilon}))u^{\epsilon} dW(s)\>\\
		&+\int_0^t\<\Delta u^{\epsilon},2g'(|u^{\epsilon}|^2) |\nabla u|^2 u^{\epsilon} dW(s)\>,\\
		II_{det} &:= \bi \Delta^2 u^{\epsilon}+\bi\lambda f(|u^{\epsilon}|^2)\Delta u^{\epsilon}dt+\bi4\lambda f'(|u^{\epsilon}|^2)Re(\bar u^{\epsilon}\nabla u^{\epsilon})\nabla u^{\epsilon}\\
		&\quad+\bi 4\lambda f''(|u^{\epsilon}|^2)(Re(\bar u^{\epsilon}\nabla u^{\epsilon}))^2u^{\epsilon}\\
		&\quad +\bi 2\lambda f'(|u^{\epsilon}|^2)Re(\bar u^{\epsilon}\Delta u^{\epsilon} )u^{\epsilon},
		\end{align*}
		and $II_{mod}$ is the summation of all terms involving the second derivative of the It\^o modified term produced by the Stratonovich integral. Here for simplicity, we omit the presentation of the explicit form for $II_{Stra}.$ 

		Taking expectation and using the Gagliardo--Nirenberg interpolation inequality $\|\nabla v\|_{L^4}\le C\|\Delta v\|^{\frac 14}\|\nabla v\|^{\frac 34}$ in $d=1$, we obtain that 
		
		\begin{align*}
		&\E\Big[ \|\Delta u^{\epsilon}(t)\|^{2}\Big]\\
		&\le \E\Big[\|\Delta u^{\epsilon}(0)\|^{2}\Big]
		+C(\lambda,p) \epsilon^{-\frac 12}\E\Big[\int_{0}^{t} \|\Delta u^{\epsilon}(r)\|(1+  \|\nabla u^{\epsilon}\|_{L^4}^{2})dr\Big]\\
		&+C(\lambda,p) \E\Big[\int_{0}^{t} \|\Delta u^{\epsilon}(r)\|\Big(\| \sum_{i} \Delta Q^{\frac 12}e_i Q^{\frac 12}e_i (g(|u^{\epsilon}|^2))^2 u^{\epsilon}\|+\| \sum_{i} |\nabla Q^{\frac 12}e_i|^2 (g(|u^{\epsilon}|^2))^2 u^{\epsilon}\|
		\\
		&+ \|\sum_{i}\nabla Q^{\frac 12}e_iQ^{\frac 12}e_ig(|u^{\epsilon}|^2)g'(|u^{\epsilon}|^2)|u^{\epsilon}|^{2}\nabla u^{\epsilon} \| +\|\sum_{i}|Q^{\frac 12}e_i|^2g(|u^{\epsilon}|^2)g'(|u^{\epsilon}|^2)|\nabla u^{\epsilon}|^2u^{\epsilon} \| \\
		&+\|\sum_{i}|Q^{\frac 12}e_i|^2(g(|u^{\epsilon}|^2)g''(|u^{\epsilon}|^2)+(g'(|u^{\epsilon}|^2))^2)|\nabla u^{\epsilon}|^2|u^{\epsilon}|^3 \|
		+\| \sum_{i} \nabla Q^{\frac 12}e_i Q^{\frac 12}e_i (g(|u^{\epsilon}|^2))^2 \nabla u^{\epsilon}\| \Big)dr\Big]\\
		&+C(\lambda,p) \E\Big[\int_{0}^{t}\sum_{i\in\mathbb N}\Big(\|g(|u^{\epsilon}|^2)\nabla u^{\epsilon}\nabla Q^{\frac 12}e_i\|^2+\|g(|u^{\epsilon}|^2)u^{\epsilon}\Delta Q^{\frac 12}e_i\|^2\\
		&+\|g'(|u^{\epsilon}|^2)|\nabla u^{\epsilon}|^2u^{\epsilon}Q^{\frac 12}e_i\|^2
		+\|g'(|u^{\epsilon}|^2)\nabla u^{\epsilon}|u^{\epsilon}|^2Q^{\frac 12}e_i\|^2\\
		&+\|g''(|u^{\epsilon}|^2)|\nabla u^{\epsilon}|^2|u^{\epsilon}|^3Q^{\frac 12}e_i\|^2+
		\|g'(|u^{\epsilon}|^2)|\nabla u^{\epsilon}|^2u^{\epsilon}Q^{\frac 12}e_i\|^2
		\Big)dr\Big]\\
		&=: \E\Big[\|\Delta u^{\epsilon}(0)\|^{2}\Big]
		+C(\lambda,p) \epsilon^{-\frac 12}\E\Big[\int_{0}^{t} \|\Delta u^{\epsilon}(r)\|(1+  \|\nabla u^{\epsilon}\|_{L^4}^{2})dr\Big]\\
		&+C(\lambda,p)  \E\Big[\int_{0}^{t} \|\Delta u^{\epsilon}(r)\| A(r) dr\Big]
		+C(\lambda,p)  \E\Big[\int_{0}^{t} B(r) dr\Big].
		\end{align*}
		Now applying the H\"older inequality, using the properties of  $g$, using the Gagliardo--Nirenberg interpolation inequality,  we obtain that for a small $\eta>0,$
		\begin{align*}
		&A(r)\le
		\sum_{i} \|\Delta Q^{\frac 12}e_i\|\| Q^{\frac 12}e_i\|_{L^{\infty}} \|(g(|u^{\epsilon}|^2))^2 u^{\epsilon}\|_{L^{\infty}}+\sum_{i} \|\nabla Q^{\frac 12}e_i\|_{L^4}^2 \|(g(|u^{\epsilon}|^2))^2 u^{\epsilon}\|_{L^{\infty}}
		\\
		&+ \sum_{i}\|\nabla Q^{\frac 12}e_i\|_{L^4}\|Q^{\frac 12}e_i\|_{L^{\infty}}\|g(|u^{\epsilon}|^2)\|_{L^{\infty}}\|\nabla u^{\epsilon} \|_{L^4} +\sum_{i}\|Q^{\frac 12}e_i\|_{L^{\infty}}^2\|\nabla u^{\epsilon}\|_{L^4}^2\|g(|u^{\epsilon}|^2)g'(|u^{\epsilon}|^2)u^{\epsilon}\|_{L^{\infty}} \\
		&+ \sum_{i}\|Q^{\frac 12}e_i\|_{L^{\infty}}^2\|\nabla u^{\epsilon}\|_{L^4}^2
		\|(g(|u^{\epsilon}|^2)g''(|u^{\epsilon}|^2)+(g'(|u^{\epsilon}|^2))^2)|u^{\epsilon}|^3 \|_{L^{\infty}}\\
		&
		+ \sum_{i} \|\nabla Q^{\frac 12}e_i\|_{L^4} \|Q^{\frac 12}e_i\|_{\infty} \|g(|u^{\epsilon}|^2)\|_{L^{\infty}}^2 \|\nabla u^{\epsilon}\|_{L^4}\\
		&\le
		\sum_{i}\Big(\|\nabla Q^{\frac 12}e_i\|_{L^4}^2+\|\nabla Q^{\frac 12}e_i\|_{L^{\infty}}^2+\|\Delta Q^{\frac 12}e_i\|^2+\|Q^{\frac 12}e_i\|_{L^{\infty}}^2\Big)
		\Big(1+\|\nabla u\|_{L^4}+\|\nabla u\|_{L^4}^2\Big)\\
		&\quad \times \Big(\|(g(|u^{\epsilon}|^2))^2 u^{\epsilon}\|_{L^{\infty}}+\|g(|u^{\epsilon}|^2)\|_{L^{\infty}} +\|g(|u^{\epsilon}|^2)g'(|u^{\epsilon}|^2)u^{\epsilon}\|_{L^{\infty}}\\
		&\quad + \|(g(|u^{\epsilon}|^2)g''(|u^{\epsilon}|^2)+(g'(|u^{\epsilon}|^2))^2)|u^{\epsilon}|^3 \|_{L^{\infty}}+\|g(|u^{\epsilon}|^2)\|_{L^{\infty}}^2 \Big)\\
		&\le  C\sum_{i}\Big(\|\nabla Q^{\frac 12}e_i\|_{L^4}^2+\|\nabla Q^{\frac 12}e_i\|_{L^{\infty}}^2+\|\Delta Q^{\frac 12}e_i\|^2+\|Q^{\frac 12}e_i\|_{L^{\infty}}^2\Big)
		\Big(1+\|\nabla u^{\epsilon}\|_{L^4}+\|\nabla u^{\epsilon}\|_{L^4}^2\Big).
		\end{align*}
		Similarly, we have that for a small $\eta>0$,
		\begin{align*}
		B(r)&\le \sum_{i}\Big(\|g(|u^{\epsilon}|^2)\|^2_{L^{\infty}}\|\nabla u^{\epsilon}\|_{L^4}^2\|\nabla Q^{\frac 12}e_i\|_{L^4}^2+\|g(|u^{\epsilon}|^2)u^{\epsilon}\|^2_{L^{\infty}}\|\Delta Q^{\frac 12}e_i\|^2\\
		&+\|g'(|u^{\epsilon}|^2) u^{\epsilon}\|_{L^{\infty}}^2\|\nabla u^{\epsilon}\|_{L^4}^4\|Q^{\frac 12}e_i\|_{L^{\infty}}^2
		+\|g'(|u^{\epsilon}|^2)|u^{\epsilon}|^2\|_{L^{\infty}}^2\|\nabla u^{\epsilon}\|_{L^4}^2\|Q^{\frac 12}e_i\|_{L^{4}}^2\\
		&+\|g''(|u^{\epsilon}|^2)|u^{\epsilon}|^3\|_{L^{\infty}}^2\|\nabla u^{\epsilon}\|_{L^4}^4\|Q^{\frac 12}e_i\|_{L^{\infty}}^2+
		\|g'(|u^{\epsilon}|^2)u^{\epsilon}\|_{L^{\infty}}^2\|\nabla u^{\epsilon}\|_{L^4}^4\|Q^{\frac 12}e_i\|_{L^{\infty}}^2
		\Big)\\
		&\le C\sum_{i}\Big(\|\nabla Q^{\frac 12}e_i\|_{L^4}^2+\|\nabla Q^{\frac 12}e_i\|_{L^{\infty}}^2+\|\Delta Q^{\frac 12}e_i\|^2+\|Q^{\frac 12}e_i\|_{L^{\infty}}^2\Big)
		\Big(1+\|\nabla u^{\epsilon}\|_{L^4}^4\Big).
		\end{align*}
		Combining the above estimates, and using the Young inequality and Gronwall inequality imply that 
		\begin{align*}
		\E\Big[\|\Delta u^{\epsilon}(t)\|^{2}\Big]\le C(u_0,T,Q,p,\eta)(1+\epsilon^{- 2}).
		\end{align*}
		Now, taking supreme over $t$, then taking expectation, and applying the Burkerholder inequality to the  $III_{sto}dW(t),$ we achieve that 
		for  sufficient small $\eta>0,$
		\begin{align*}
		\E\Big[\sup_{t\in[0,\tau]}\|\Delta u^{\epsilon}(t)\|^{2}\Big]\le C(u_0,T,Q,\eta)(1+\epsilon^{-2}).
		\end{align*}
	\qed
	
	\textbf{[Proof of Proposition \ref{wei-hol-2}]}
		We follow the steps in the proof of Proposition \ref{wei-hol-1} to present the proof in the case of $p=2$.  For convenience, we present the proof for the multiplicative noise case. 
		Applying the It\^o formula to $\|u^{\epsilon}\|_{L^2_{\alpha}}^{2}=\int_{\mathbb R^d}(1+|x|^2)^{\alpha}|u^{\epsilon}|^2dx$,
		using integration by parts, then taking supreme over $t$, and applying Burkerholder inequality, we deduce that 
		\begin{align*}
		\E\Big[\sup_{t\in [0,T]}\|u^{\epsilon}(t)\|_{L^2_{\alpha}}^{2}\Big]
		&\le \E\Big[\|u_0\|_{L^2_{\alpha}}^{2}\Big]
		+2\alpha \E \Big[\int_0^T \Big| \<(1+|x|^2)^{\alpha-1} x u^{\epsilon}(s),\nabla u^{\epsilon}\> \Big|ds\Big]\\
		&+\E \Big[\sup_{t\in [0,T]}\Big|\int_0^t\<(1+|x|^2)^{\alpha}u^{\epsilon}(s),\mathbf i g(|u^{\epsilon}(s)|^2)u^{\epsilon}(s)dW(s)\>\Big|\Big]\\
		&\le \E\Big[\|u_0\|_{L^2_{\alpha}}^{2}\Big]
		+C_{\alpha} \E \Big[\int_0^T\Big| \<(1+|x|^2)^{\alpha-1} x u^{\epsilon}(s),\nabla u^{\epsilon}\> \Big|ds\Big]\\
		&+C \E \Big[\int_0^T\sum_{i\in\mathbb N^+}\Big|\<(1+|x|^2)^{\alpha}u^{\epsilon}(s), \mathbf ig(|u^{\epsilon}(s)|^2)u^{\epsilon}(s) Q^{\frac 12}e_i\>\Big|^2\Big|ds\Big].
		\end{align*}
		By H\"older's inequality, for $\alpha\in (1,2]$, we have that 
		\begin{align*}
		\Big| \<(1+|x|^2)^{\alpha-1} x u^{\epsilon},\nabla u^{\epsilon}\> \Big|
		&\le C\|u^{\epsilon}\|_{L^{2}_{\alpha}}\|(1+|x|^2)^{\frac \alpha 2-\frac 12}\nabla u^{\epsilon}\|.
		\end{align*}
		Integration by parts and H\"older's inequality yield that for some small $\eta>0,$
		\begin{align*}
		&\|(1+|x|^2)^{\frac \alpha 2-\frac 12}\nabla u^{\epsilon}\|^2=\<(1+|x|^2)^{\alpha -1} \nabla u^{\epsilon}, \nabla u^{\epsilon}\>\\
		&=-\<(1+|x|^2)^{\alpha -1}u^{\epsilon}, \Delta u^{\epsilon}\>
		-2(\alpha-1)\<(1+|x|^2)^{\alpha -2} x \nabla u^{\epsilon},  u^{\epsilon}\>\\
		&\le \|u^{\epsilon}\|_{L^2_{\max(2\alpha-2,0)}} \|\Delta u^{\epsilon}\|+ C(\eta)|\alpha-1|\|u^{\epsilon}\|^2_{L^2_{\alpha}}+\eta |\alpha-1|\|\nabla u^{\epsilon}\|^2_{L^2_{\max(\alpha-3,0)}}.
		\end{align*}
		Combining the above estimates in Proposition \ref{wei-hol-1} and using Young's inequality, we achieve that
		\begin{align*}
		\E\Big[\sup_{t\in[0,T]}\|u^{\epsilon}(t)\|_{L^2_{\alpha}}^{2}\Big]&\le  e^{CT}(1+\epsilon^{-1}),\; \text{if}\; \alpha \in (1,\frac 32],\\
		\E\Big[\sup_{t\in[0,T]}\|u^{\epsilon}(t)\|_{L^2_{\alpha}}^{2}\Big]&\le e^{CT}(1+\epsilon^{-\frac 32}),\; \text{if}\; \alpha \in [\frac 32,2),\\
		\E\Big[\sup_{t\in[0,T]}\|u^{\epsilon}(t)\|_{L^2_{\alpha}}^{2}\Big]&\le e^{CT}(1+\epsilon^{-2}),\; \text{if}\; \alpha =2.\\
		\end{align*} 	
	\qed

	\bibliographystyle{plain}
	\bibliography{bib}
\end{document}